\documentclass[11pt]{article}

\usepackage{oldlfont}
\usepackage{amsfonts}
\usepackage{graphicx}
\usepackage{eepic}

\newcommand{\be}{\begin{equation}}
\newcommand{\ee}{\end{equation}}
\newcommand{\bea}{\begin{eqnarray}}
\newcommand{\eea}{\end{eqnarray}}
\newcommand{\beas}{\begin{eqnarray*}}
\newcommand{\eeas}{\end{eqnarray*}}
\newcommand{\ba}{\begin{array}}
\newcommand{\ea}{\end{array}}

\newcommand{\<}  {\langle}
\renewcommand{\>}{\rangle}
\newcommand{\field}[1]{\mathbb{#1}}

\newcommand{\diam}{{\mathrm diam}}
\newcommand{\grad}{\nabla}
\renewcommand{\div}{{\mathrm div}}
\newcommand{\curl}{{\mathrm curl}}

\newcommand{\G}{\Gamma}
\newcommand{\g}{\gamma}
\newcommand{\eps}{\varepsilon}

\newcommand{\bolda}{{\mathbf a}}

\newcommand{\bn}{{\mathbf n}}
\newcommand{\bq}{{\mathbf q}}
\newcommand{\bu}{{\mathbf u}}
\newcommand{\bv}{{\mathbf v}}
\newcommand{\bw}{{\mathbf w}}
\newcommand{\bx}{{\mathbf x}}

\newcommand{\boldb}{{\mathbf b}}
\newcommand{\boldf}{{\mathbf f}}
\newcommand{\bA}{{\mathbf A}}
\newcommand{\bB}{{\mathbf B}}

\newcommand{\bV}{{\mathbf V}}
\newcommand{\bW}{{\mathbf W}}
\newcommand{\bX}{{\mathbf X}}
\newcommand{\bH}{{\mathbf H}}
\newcommand{\bL}{{\mathbf L}}
\newcommand{\bcurl}{{\mathbf curl}}
\newcommand{\bzero}{{\mathbf 0}}

\newcommand{\gradg}{\nabla_{\G}}
\newcommand{\divg}{{\mathrm div}_{\G}}
\newcommand{\curlg}{{\mathrm curl}_{\G}}
\newcommand{\bcurlg}{{\mathbf curl}_{\G}}

\newcommand{\bnu}{\hbox{\mathversion{bold}$\nu$}}
\newcommand{\bxi}{\hbox{\mathversion{bold}$\xi$}}
\newcommand{\bPsi}{\hbox{\mathversion{bold}$\Psi$}}

\newcommand{\bphi}{\hbox{\mathversion{bold}$\phi$}}
\newcommand{\btau}{\hbox{\mathversion{bold}$\tau$}}

\newcommand{\CH}{{\cal H}}
\newcommand{\CI}{{\cal I}}
\newcommand{\CL}{{\cal L}}
\newcommand{\CM}{{\cal M}}
\newcommand{\CP}{{\cal P}}
\newcommand{\bCP}{\hbox{\mathversion{bold}$\cal P$}}

\newcommand{\CT}{{\cal T}}

\newcommand{\stack}[2]{\mathrel{\mathop{#2}\limits_{\scriptstyle #1}}}
\newcommand{\tH}{\tilde H}

\newcommand{\ttu}{\widetilde{\widetilde{u^\circ}}}
\newcommand{\ttv}{\widetilde{\widetilde{v^\circ\,}}}

\newtheorem{theorem}{Theorem} [section]
\newtheorem{lemma}{Lemma} [section]

\newtheorem{prop}{Proposition} [section]
\newtheorem{remark}{Remark} [section]
\newenvironment{proof}{\noindent\textbf{Proof.}\ }
              {\nopagebreak\hbox{ }\hfill$\Box$\bigskip}

\title{On the convergence of the $hp$-BEM with quasi-uniform meshes
       for the electric field integral equation on polyhedral surfaces
\thanks{Supported by EPSRC under grant no. EP/E058094/1.}}

\author{Alexei Bespalov
\thanks{Department of Mathematical Sciences, Brunel University,
        Uxbridge, West London UB8 3PH, UK.
        Email: {\tt albespalov@yahoo.com}\
        Supported by the INTAS Young Scientist Fellowship grant
        (project no. 06-1000014-5945).}
        \and
        Norbert Heuer
\thanks{Facultad de Matem\'aticas, Pontificia Universidad Cat\'olica de Chile,
        Avenida Vicu\~na Mackenna 4860, Santiago, Chile.
        Email: {\tt nheuer@mat.puc.cl}}
        }
\begin{document}
\date{}
\maketitle

\begin{abstract}
In this paper the $hp$-version of the boundary element method is applied to the electric
field integral equation on a piecewise plane (open or closed) Lipschitz surface.
The underlying meshes are supposed to be quasi-uniform.
We use $\bH(\div)$-conforming discretisations with quadrilateral elements of
Raviart-Thomas type and establish quasi-optimal convergence of $hp$-approximations.
Main ingredient of our analysis is a new $\tilde\bH^{-1/2}(\div)$-conforming
$p$-interpolation operator that assumes only $\bH^r \cap \tilde\bH^{-1/2}(\div)$-regularity
($r > 0$) and for which we show quasi-stability with respect to polynomial degrees.

\bigskip
\noindent
{\em Key words}: $hp$-version with quasi-uniform meshes, electric field integral equation,
                 time-harmonic electro-magnetic scattering,
                 boundary element method

\noindent
{\em AMS Subject Classification}: 65N38, 65N12, 78M15, 41A10
\end{abstract}

\section{Introduction and formulation of the problem} \label{sec_intro}
\setcounter{equation}{0}

In this paper we prove convergence of the $hp$-version of the boundary element
method (BEM) for the electric field integral equation (EFIE) on piecewise
plane (open or closed) surfaces discretised by quasi-uniform meshes.
The EFIE is a boundary integral equation that
represents a boundary value problem for the time-harmonic Maxwell's equations in the
exterior domain. It models the scattering of electro-magnetic waves
at a perfectly conducting body (the scatterer). The solution to the EFIE is
the induced electric surface current (a tangential vector field) on the surface
of the scatterer, see, e.g., \cite{Nedelec_01_AEE}.
If the scatterer is a thin object (i.e., its thickness is small in
comparison to the wave length), then it can be modelled as an open
surface (a sub-manifold with boundary) in ${\field{R}}^3$. Our analysis covers this
theoretically challenging case, which has important applications
(e.g., antenna problems).

The basis of our BEM is a variational formulation of the EFIE,
called Rumsey's principle. For smooth surfaces, its boundary element
discretisation has been studied by Bendali in 1984, see \cite{Bendali_84_NA1,Bendali_84_NA2}.
Progress in the numerical analysis of the EFIE on Lipschitz surfaces
has been achieved relatively recently and was inspired by the study of traces
of functional spaces that govern Maxwell's equations in Lipschitz domains
\cite{BuffaCS_02_THL}. The main challenges in this analysis concern the solvability
and quasi-optimal convergence of approximations and a priori error estimation
in the energy norm. In the framework of the $h$-version of the BEM, i.e.,
for discretisations with elements of fixed order on refined meshes, these issues
were addressed in \cite{BuffaCS_02_BEM,HiptmairS_02_NBE,BuffaHvPS_03_BEM}
(for polyhedral surfaces) and in \cite{BuffaC_03_EFI} (for open Lipschitz
surfaces). We note that in \cite{HiptmairS_02_NBE,BuffaHvPS_03_BEM,BuffaC_03_EFI}
the authors focused on conforming discretisations of Rumsey's principle with
Raviart-Thomas (or Brezzi-Douglas-Marini) boundary elements, whereas 
in \cite{BuffaCS_02_BEM} a mixed formulation
utilising standard (continuous) basis functions was used. In this paper we follow
the former approach, called the natural boundary element method for the EFIE.

While in the $p$-version of the BEM the mesh is fixed and approximations
are improved by increasing polynomial degrees, the $hp$-version combines both
mesh refinement and increase of polynomial degrees.
In our previous paper \cite{BespalovH_NpB} we analysed the natural $p$-BEM
for the EFIE on a plane open surface with polygonal boundary. We have proved
convergence of the $p$-version with Raviart-Thomas (RT) parallelogram elements
and derived an a priori error estimate which takes into account the strong
singular behaviour of the solution at edges and corners of the surface.
With the present note we prove the convergence of the natural $hp$-BEM
on polyhedral and piecewise plane open surfaces discretised by
quasi-uniform meshes of quadrilateral (in general, curvilinear) elements.
We emphasize that RT-spaces are used for both affine and non-affine quadrilateral
meshes.

In order to prove convergence of approximations for the EFIE, one usually relies
on properties of the continuous and discrete Helmholtz decompositions, and on the
proximity in some sense of the discrete decompositions to the continuous one,
see \cite{BuffaC_03_EFI,BuffaHvPS_03_BEM}. The key property is the orthogonality
of decompositions, and the main tool in the analysis is an appropriate interpolation
operator (a projector onto the corresponding polynomial space).
In \cite{BuffaC_03_EFI, BuffaHvPS_03_BEM}, $\bL^2$-orthogonal discrete
decompositions mimicking the Helmholtz decomposition of the energy space were
analysed for finite dimensional subspaces based on Raviart-Thomas and
Brezzi-Douglas-Marini (BDM) boundary elements.
It has been proved that these discrete decompositions are sufficiently close
to the continuous one as the mesh parameter $h$ tends to zero (i.e., for
the $h$-version of the BEM). The main tools in the proofs were
the standard $\bH(\div)$-conforming RT or BDM interpolation operators.
However, it turns out that a generalisation of that approach to the $p$- and
the $hp$-versions is not straightforward when sticking to both the $\bL^2$-orthogonality
of decompositions and classical interpolation operators.
In \cite{BespalovH_NpB}, for the $p$-version, we employed an $\tilde\bH^{-1/2}$-orthogonality
of the Helmholtz decomposition, while using the classical RT interpolation
operator. Unfortunately, the extension of this approach to polyhedral surfaces
and even to piecewise plane screens does not seem to be easy. In particular, the low
regularity of the Laplace-Beltrami operator on polyhedral surfaces is not enough
to prove stability (with respect to polynomial degrees) of the classical
Raviart-Thomas interpolation operator. That is why in
this paper we use an alternative approach: we adhere to the $\bL^2$-orthogonality
of decompositions but utilise a non-classical $\tilde\bH^{-1/2}(\div)$-conforming
interpolation operator (for RT-elements). Our construction of this operator (see
Section~\ref{sec_inter}) is much in the spirit of \cite{DemkowiczB_03_pIE},
where $H^1$- and $\bH(\curl)$-conforming projection-based interpolation
operators were introduced and analysed. We, however, employ the $\tH^{-1/2}$-projection
for the divergence term, and thus need the corresponding inner product to be written
in an appropriate explicit form. Moreover, the use of an appropriate scaling argument
allows to prove the convergence result in the framework of the $hp$-BEM.

We will denote by $\G$ a piecewise plane (open or closed) Lipschitz surface in ${\field{R}}^3$.
Let us introduce Rumsey's formulation of the electric field integral equation on $\G$.
For a given wave number $k>0$ and a scalar function $v$ (resp., tangential vector field $\bv$)
we define the single layer operator $\Psi_k$ (resp., $\bPsi_k$) by
\beas
     & \displaystyle{
     \Psi_k v(x) =
     \frac 1{4\pi}\int_\G v(y) \frac {e^{ik|x-y|}}{|x-y|}\,dS_y,\qquad
     x \in {\field{R}}^3 \backslash \G}
     &
     \\[5pt]
     & \displaystyle{
     \Big(\hbox{resp.,}\quad
     \bPsi_k \bv(x) =
     \frac 1{4\pi}\int_\G \bv(y) \frac {e^{ik|x-y|}}{|x-y|}\,dS_y,\qquad
     x \in {\field{R}}^3 \backslash \G\Big).}
     &
\eeas
Let $\bL^2_t(\G)$ be the space of two-dimensional, tangential, square integrable
vector fields on $\G$. By $\gradg$ (resp., $\divg$) we denote the surface gradient
(resp., surface divergence) acting on scalar functions
(resp., tangential vector fields) on $\G$. We will need the following space:
\[
  \bX = \bH^{-1/2}(\divg,\G) :=
        \{\bu \in \bH^{-1/2}_{\|}(\G);\; \divg\,\bu \in H^{-1/2}(\G)\}
\]
if $\G$ is a closed surface, and
\beas
     \bX = \tilde\bH^{-1/2}_0(\divg,\G)
     & := &
     \{\bu \in \tilde\bH^{-1/2}_{\|}(\G);\ \divg\,\bu \in \tH^{-1/2}(\G)\ \hbox{and}
     \\[3pt]
     &   &
     \quad
     \<\bu,\gradg v\> + \<\divg\,\bu,v\> = 0\quad
     \hbox{for all}\ \ v \in C^{\infty}(\G)\}
\eeas
if $\G$ is an open surface.
In the latter definition the brackets $\<\cdot,\cdot\>$ denote
dualities associated with $\bH^{1/2}_{\|}(\G)$ and $H^{1/2}(\G)$, respectively.
For definitions of the space $C^{\infty}(\G)$ and the Sobolev spaces on $\G$
see \S\ref{sec_spaces}.

Throughout, we use boldface symbols for vector fields.
The spaces (or sets) of vector fields are also denoted in boldface
(e.g., $\bH^s(\G) = (H^s(\G))^3$), with their norms and inner
products being defined in~\S\ref{sec_spaces}.

Let $\bX'$ be the dual space of $\bX$ (with $\bL^2_t(\G)$ as pivot space).
Now, for a given tangential vector field $\boldf\in \bX'$
($\boldf$ represents the excitation by an incident wave), Rumsey's formulation
reads as:
{\em find a complex tangential field $\bu\in\bX$ such that}
\be \label{bie_var}
    a(\bu,\bv)
    :=
    \<\g_{\rm tr}(\Psi_k\divg\,\bu),\divg\,\bv\> - k^2 \<\pi_{\tau}(\bPsi_k\bu), \bv\>
    = \<\boldf, \bv\> \quad\forall \bv \in \bX.
\ee
Here $\g_{\rm tr}$ is the standard trace operator, and $\pi_{\tau}$ denotes the
tangential components trace mapping (see \S\ref{sec_spaces} for definitions and
properties of these operators).

The rest of the paper is organised as follows.
In the next section we define the $hp$-version of the BEM for the EFIE and formulate
the main result (Theorem~\ref{thm_solve}), which states the unique solvability and
quasi-optimal convergence of this approximation method. Section~\ref{sec_prelim}
gives necessary preliminaries: in \S\ref{sec_spaces} we recall definitions of functional
spaces of scalar functions and vector fields;
then in \S\ref{sec_norms} we introduce some equivalent norms
in the Sobolev spaces $H^s$ and $\tH^s$ ($s = \pm\frac 12$) on the reference element
and derive expressions for corresponding inner products; some auxiliary lemmas
are collected in \S\ref{sec_aux}. In Section~\ref{sec_dec} we discuss the continuous
Helmholtz decomposition of the energy space $\bX$ and its discrete counterpart.
Section~\ref{sec_inter} is devoted to interpolation operators.
First, we recall some known operators such as
the $L^2$- and $\tH^{-1/2}$-projectors,
the $H^1$- and $\bH(\curl)$-conforming projection-based interpolation operators.
Then we introduce an $\tilde\bH^{-1/2}(\div)$-conforming interpolation operator
and study its properties. In particular, we prove the quasi-stability of this operator
and its commutativity with the $\tH^{-1/2}$-projector. The $\tilde\bH^{-1/2}(\div)$-conforming
interpolation operator and $\bL^2$-orthogonal discrete Helmholtz decompositions
are the main tools in the proof of Theorem~\ref{thm_solve}, which is given in
Section~\ref{sec_proof_solve}. We conclude the paper with Section~\ref{sec_conclusion},
where we comment on extensions and open problems.

Throughout the paper, $C$ denotes a generic positive constant which is independent
of $h$, $p$, and involved functions, unless stated otherwise.

\section{The $hp$-version of the BEM and the main result} \label{sec_p-BEM}
\setcounter{equation}{0}

For the approximate solution of (\ref{bie_var}) we apply the $hp$-version of the BEM based
on Galerkin discretisations with Raviart-Thomas spaces on quasi-uniform meshes.
In what follows, $h > 0$ and $p \ge 1$ will always specify the mesh parameter and
a polynomial degree, respectively.
For any $\Omega\subset {\field{R}}^n$ we will denote
$\rho_{\Omega} = \sup \{\diam (B);\; \hbox{$B$ is a ball in $\Omega$}\}$.
Furthermore, we will denote by $K=(0,1)^2$ the reference square.
The sides of $K$ will be denoted by $\ell_i$ ($i = 1,\ldots,4$).

Let $\CT = \{\Delta_h\}$ be a family of meshes
$\Delta_h = \{\G_j;\; j=1,\ldots,J\}$ on $\G$, where the elements $\G_j$ are open
quadrilaterals (in general, curvilinear ones) satisfying the following standard
assumptions (below we denote $h_j=\diam (\G_j)$ for any $\G_j\in \Delta_h$):

\begin{enumerate}

\item[i)]
$\bar\G = \cup_{j=1}^{J} \bar\G_j$;
the intersection of any two quadrilaterals $\bar\G_j,\,\bar\G_k$ ($j \not= k$)
is either a common vertex, an entire side, or empty;

\item[ii)]
The elements are shape regular, i.e., there exists a positive constant $C$
independent of $h = \max\limits_j h_j$ such that for
any $\G_j\in\Delta_h$ and arbitrary $\Delta_h\in \CT$ there holds
$h_j \le C\,\rho_{\G_j}$.

\item[iii)]
Any element $\G_j$ is the image of the reference square $K$, more precisely
\[
  \bar\G_j = T_j(\bar K),\quad
  \bx = T_j(\bxi),\ 
  \bx = (x_1,x_2) \in \bar\G_j,\ \bxi = (\xi_1,\xi_2) \in \bar K,
\]
where $T_j$ is sufficiently smooth one-to-one mapping with sufficiently
smooth inverse $T_j^{-1}:\,\G_j \rightarrow K$.
The Jacobian matrix of $T_j$ is denoted by $DT_j(\bxi)$, it is supposed to be
invertible for any $\bxi \in K$, and $DT_j^{-1}(\bx) = (DT_j(\bxi))^{-1}$.
We assume that
\be \label{Jacobian}
    |\hbox{det}(DT_j(\bxi))| \simeq h_j^2\qquad \hbox{for any $\bxi \in K$},
\ee
and there exist positive constants $C$ independent of $h$ such that
for $k = 1,2$ there holds
\be \label{mesh-reg}
    \displaystyle{
    \sup_{\bxi \in K} \|D^k T_j(\bxi)\|_{\CL_k({\field{R}}^2,{\field{R}}^2)}
    }
    \le C\,h_j^k,\qquad
    \displaystyle{
    \sup_{\bx \in \G_j} \|D^k T_j^{-1}(\bx)\|_{\CL_k({\field{R}}^2,{\field{R}}^2)}
    }
    \le C\,h_j^{-1}.
\ee
Here $D^k T_j(\bxi)$ (resp., $D^k T_j^{-1}(\bx)$) denotes the $k$-th (Fr{\'e}chet) derivative
of $T_j$ (resp., $T_j^{-1}$) and $\|\cdot\|_{\CL_k(X,Y)}$ is the operator
norm in the space $\CL_k(X,Y)$ of continuous $k$-linear mappings from $X^k$ into $Y$.

\item[iv)]
if $\bar\G_j \cap \bar\G_k$ is an entire side $\ell$,
then denoting by $T^{\ell}_j$ (resp., $T^{\ell}_k$) the restriction
of $T_j$ (resp., $T_k$) to the corresponding side 
$T_j^{-1}(\ell)$ (resp., $T_k^{-1}(\ell)$) of the reference element $K$, one has 
$T^{\ell}_j \equiv T^{\ell}_k$ as mappings of the unit interval $(0,1)$ onto $\ell$.

\end{enumerate}

\begin{remark} \label{rem_mesh-reg}
Assumptions {\rm (\ref{Jacobian})}, {\rm (\ref{mesh-reg})} above are always satisfied
for affine families of elements (i.e., for meshes of parallelograms).
For curvilinear elements these assumptions
are satisfied, for instance, if the elements tend to be affine as
$h \rightarrow 0$ (see {\rm \cite[Section~4.3]{Ciarlet_78_FEM}}).
In this case {\rm (\ref{Jacobian})}, {\rm (\ref{mesh-reg})} follow
from assumption ii) and the smoothness of $T_j$, provided that
$h_j$ is small enough.
\end{remark}

In this paper we consider a family $\CT$ of quasi-uniform meshes $\Delta_h$
on $\G$ in the sense that there exists a positive constant
$C$ independent of $h$ such that for any $\G_j\in\Delta_h$ and
arbitrary $\Delta_h\in \CT$ there holds $h \le C\,h_j$.

The mapping $T_j$ introduced above is used to associate the scalar function $u$ defined
on the real element $\G_j$ with the function $\hat u$ defined on the reference
element $K$:
\[
  u = \hat u \circ T_j^{-1}\ \ \hbox{on $\G_j$\qquad and}\qquad
  \hat u = u \circ T_j\ \ \hbox{on $K$}.
\]
Any vector-valued function $\hat\bv$ defined on $K$ is transformed to the function
$\bv$ on $\G_j$ by using the Piola transformation:
\be \label{Piola}
    \bv = \CM_j(\hat \bv) = \hbox{$\frac {1}{J_j}$} DT_j \hat\bv \circ T_j^{-1},
    \quad
    \hat\bv = \CM_j^{-1}(\bv) = J_j DT_j^{-1} \bv \circ T_j,
\ee
where $J_j$ is the determinant $J_j(\bxi) := \hbox{det}(DT_j(\bxi))$.

Let us introduce the needed polynomial sets.
By $\CP_p(I)$ we denote the set of polynomials of degree $\le p$ on
an interval $I\subset\field{R}$, and $\CP^0_p(I)$ denotes the subset of $\CP_p(I)$
which consists of polynomials vanishing at the end points of $I$.
In particular, these two sets will be used for edges $\ell_i \subset \partial K$.

Further, $\CP_{p_1,p_2}(K)$ denotes the set of polynomials on $K$ of degree
$\le p_1$ in $\xi_1$ and degree $\le p_2$ in $\xi_2$. For $p_1 = p_2 = p$ we denote
$\CP_{p}(K) = \CP_{p,p}(K)$. The corresponding set of polynomial (scalar) bubble
functions on $K$ is denoted by $\CP^0_{p}(K)$.

Let us denote by $\bCP^{\rm RT}_p(K)$ the RT-space of order $p\ge 1$ on $K$
(see, e.g., \cite{BrezziF_91_MHF, RobertsT_91_MHM}), i.e.,
\[
  \bCP^{\rm RT}_p(K) = \CP_{p,p-1}(K) \times \CP_{p-1,p}(K).
\]
The subset of $\bCP^{\rm RT}_p(K)$ which consists of vector-valued polynomials
with vanishing normal trace on the boundary $\partial K$
(vector bubble-functions) will be denoted by $\bCP^{\rm RT,0}_p(K)$.

Then using transformations (\ref{Piola}), we set
\be \label{Xp}
    \bX_{hp} := \{\bv \in \bX^0;\;
                      \CM_j^{-1}(\bv|_{\G_j}) \in \bCP^{\rm RT}_p(K),\ j=1,\ldots,J\},
\ee
where the space $\bX^0 \subset \bX$ is defined in \S\ref{sec_spaces}
($\bX^0 = \bH(\divg,\G)$ if $\G$ is closed and
$\bX^0 = \bH_0(\divg,\G)$ if $\G$ is an open surface).
We will denote by $N = N(h,p)$ the dimension of the discrete space $\bX_{hp}$.
One has $N \simeq h^{-2}$ for fixed $p$ and $N \simeq p^2$ for fixed $h$.

The $hp$-version of the Galerkin BEM for the EFIE reads as:
{\em Find $\bu_{hp}\in \bX_{hp}$ such that}
\be \label{BEM}
    a(\bu_{hp},\bv) = \<\boldf, \bv\> \quad \forall \bv \in \bX_{hp}.
\ee

Let us formulate the result which states the unique solvability of (\ref{BEM})
and quasi-optimal convergence of the $hp$-version of the BEM for the EFIE.

\begin{theorem} \label{thm_solve}
There exists $N_0 \ge 1$ such that for any $\boldf \in \bX'$ and for arbitrary
mesh-degree combination satisfying $N(h,p) \ge N_0$ the discrete problem {\rm (\ref{BEM})}
is uniquely solvable and the $hp$-version of the Galerkin BEM generated by RT-elements
converges quasi-optimally, i.e.,
\[
  \|\bu - \bu_{hp}\|_{\bX} \le
  C \inf\{\|\bu - \bv\|_{\bX};\; \bv\in \bX_{hp}\}.
\]
Here, $\bu \in \bX$ is the solution of {\rm (\ref{bie_var})},
$\bu_{hp} \in \bX_{hp}$ is the solution of {\rm (\ref{BEM})},
$\|\cdot\|_{\bX}$ denotes the norm in $\bX$, and
$C>0$ is a constant independent of $h$ and $p$.
\end{theorem}

The proof of Theorem~\ref{thm_solve} is given in Section~\ref{sec_proof_solve} below.

\section{Preliminaries} \label{sec_prelim}
\setcounter{equation}{0}

\subsection{Functional spaces, norms, and inner products} \label{sec_spaces}

First, let us recall the Sobolev spaces and norms for scalar functions
on a Lipschitz domain $\Omega\subset{\field{R}}^n$, see \cite{LionsMagenes}.
To that end we will need the space $C^\infty(\Omega)$ of infinitely
differentiable functions in $\Omega$ and its subspace
$C_0^\infty(\Omega) \subset C^\infty(\Omega)$ which consists of functions with compact
support in $\Omega$.

For an integer $s$, let $H^s(\Omega)$
be the closure of $C^\infty(\Omega)$ with respect to the norm
\[
   \|u\|_{H^s(\Omega)}^2 = \|u\|_{H^{s-1}(\Omega)}^2 + |u|_{H^s(\Omega)}^2
   \quad (s\ge 1),
\]
where
\[
   |u|_{H^s(\Omega)}^2 = \int_{\Omega}|D^su(x)|^2\,dx,\quad\mbox{and}\quad
   H^0(\Omega) = L^2(\Omega).
\]
Here, $|D^su(x)|^2=\sum_{|\alpha|=s}|D^\alpha u(x)|^2$ in the usual notation with
multi-index $\alpha=(\alpha_1,\ldots,\alpha_n)$ and with respect to Cartesian
coordinates $x=(x_1,\ldots,x_n)$.
For a positive non-integer $s=m+\sigma$ with integer $m\ge 0$ and
$0<\sigma<1$, the norm in $H^s(\Omega)$ is
\[
   \|u\|_{H^s(\Omega)}^2 = \|u\|_{H^m(\Omega)}^2 + |u|_{H^s(\Omega)}^2
\]
with semi-norm
\[
   |u|_{H^s(\Omega)}^2
   =
   \sum_{|\alpha|=m}\int_{\Omega}\int_{\Omega}
   \frac {|D^\alpha u(x) - D^\alpha u(y)|^2}{|x-y|^{n+2\sigma}}\,dx\,dy.
\]
The Sobolev spaces $\tH^s(\Omega)$ for $s \in (0,1)$ and for a bounded
Lipschitz domain $\Omega$ are defined by interpolation.
We use the real K-method of interpolation (see \cite{LionsMagenes}) to define
\[
   \tilde H^{s}(\Omega) = \Big(L^2(\Omega), H_0^t(\Omega)\Big)_{\frac st,2}
   \quad (1/2 < t \le 1,\ 0<s<t).
\]
Here, $H_0^t(\Omega)$  ($0<t\le 1$) is the completion of $C_0^\infty(\Omega)$ in
$H^t(\Omega)$ and we identify $H_0^1(\Omega)$ and $\tilde H^1(\Omega)$.
Note that the Sobolev spaces $H^s(\Omega)$ also satisfy the interpolation property
\[
   H^s(\Omega) = \Big(L^2(\Omega), H^1(\Omega)\Big)_{s,2}\quad (0<s<1)
\]
with equivalent norms. Furthermore, the semi-norm $|\cdot|_{H^1(\Omega)}$
is a norm in $\tH^1(\Omega)$ due to the Poincar{\' e} inequality.
For the $L^2(\Omega)$-norm we will use the notation $\|\cdot\|_{0,\Omega}$.

For $s\in[-1,0)$ the Sobolev spaces and their norms are defined by duality
with $L^2(\Omega) = H^0(\Omega) = \tH^0(\Omega)$ as pivot space:
\[
   H^s(\Omega) = \big(\tilde H^{-s}(\Omega)\big)',\quad
   \tilde H^s(\Omega) = \big(H^{-s}(\Omega)\big)',
\]
\be \label{dual_norm}
    \|u\|_{H^{s}(\Omega)} = \sup_{0\not=v\in\tH^{-s}(\Omega)}
    {|\<u,v\>| \over{\|v\|_{\tH^{-s}(\Omega)}}},
    \quad
    \|u\|_{\tH^{s}(\Omega)} = \sup_{0\not=v\in H^{-s}(\Omega)}
    {|\<u,v\>| \over{\|v\|_{H^{-s}(\Omega)}}},
\ee
where
\[
  \<u,v\> = \<u,v\>_{0,\Omega} := \int_\Omega u(x) \bar v(x) dx
\]
denotes the extension of the $L^2(\Omega)$-inner product by duality
(and $\bar v$ is the complex conjugate of $v$).

Now, let $\G$ be a piecewise smooth (open or closed) Lipschitz surface
in ${\field{R}}^3$. We will assume that $\G$ has plane faces $\G^{(i)}$
($i = 1,\ldots,\CI$; without loss of generality it is assumed that $\CI > 1$)
and straight edges $e_{ij} = \bar\G^{(i)} \cap \bar\G^{(j)} \not= \mbox{\o}$ ($i \not= j$).
If $\G$ is a closed surface, we will denote by
$\Omega$ the Lipschitz polyhedron bounded by $\G$, i.e., $\G=\partial\Omega$.
In the case of an open surface $\G$, we first introduce a piecewise plane
closed Lipschitz surface $\tilde\G$ which contains $\G$, and then denote by
$\Omega$ the Lipschitz polyhedron bounded by $\tilde\G$, i.e., $\tilde\G=\partial\Omega$.
For each face $\G^{(i)} \subset \G$ there exists a constant unit normal vector $\bnu_i$,
which is an outer normal vector to $\Omega$.
These vectors are then blended into a unit normal vector
$\bnu$ defined almost everywhere on~$\G$.
For each pair of indices $i,j =  1,\ldots,\CI$ such that
$\bar\G^{(i)} \cap \bar\G^{(j)} = e_{ij}$
we consider unit vectors $\btau_{ij}$, $\btau_i^{(j)}$, and $\btau_j^{(i)}$
such that $\btau_{ij} \| e_{ij}$, $\btau_{i}^{(j)} = \btau_{ij} \times \bnu_i$,
and $\btau_{j}^{(i)} = \btau_{ij} \times \bnu_j$.
Since each $\G^{(i)}$ can be identified with a bounded subset in ${\field{R}}^2$,
the pair $(\btau_{i}^{(j)},\btau_{ij})$ is an orthonormal basis of the plane generated
by $\G^{(i)}$.

Let $\G$ be a closed surface. Then $\G=\partial\Omega$ is
locally the graph of a Lipschitz function. Since the Sobolev
spaces $H^s$ for $|s| \le 1$ are invariant under Lipschitz (i.e., $C^{0,1}$)
coordinate transformations, the spaces $H^s(\G)$ with $|s| \le 1$
are defined in the usual way via a partition of unity subordinate to a finite family
of local coordinate patches (see \cite{McLean_00_SES}).
Due to this definition, the properties of Sobolev spaces on Lipschitz domains
in ${\field{R}}^n$ carry over to Sobolev spaces on Lipschitz surfaces.
If $\G$ is an open surface, then the Sobolev spaces
$H^s(\G)$, $\tH^s(\G)$ for $|s| \le 1$ and $H^s_0(\G)$ for $0 < s \le 1$
are constructed in terms of the Sobolev spaces $H^s(\tilde\G)$ on
a closed Lipschitz surface $\tilde\G \supset \G$ (see \cite{McLean_00_SES}). 
Note that the spaces $H^s(\G^{(i)})$ and $\tH^s(\G^{(i)})$ on each face
$\G^{(i)}$ are well-defined for any $s \ge -1$.

For $s > 1$ we define the space $H^s(\G)$ in the following piecewise fashion
(hereafter, $u_i$ denotes the restriction of $u$ to the face $\G^{(i)}$):
\[
  H^s(\G) := \{u \in H^1(\G);\; u_i \in H^s(\G^{(i)}),\ 
               i = 1,\ldots,\CI\}.
\]
This space is equipped with its natural norm
\[
  \|u\|_{H^s(\G)} := 
  \bigg(\|u\|^2_{H^1(\G)} +
  \sum\limits_{i=1}^{\CI} \|u_i\|^2_{H^{s}(\G^{(i)})}\bigg)^{\frac 12}.
\]
Besides the above, we will need the following spaces:
\[
  H^s_{*}(\G) := \{u \in H^s(\G);\; \<u,1\>_{0,\G} = 0\},
\]
where $s \ge -1$ if $\G$ is closed, and $s > -\frac 12$ if $\G$ is an open surface.

We will denote by $\g_{\rm tr}$ the standard trace operator,
$\g_{\rm tr}(u) = u|_{\G}$, $u \in C^\infty(\bar\Omega)$.
For $s \in (0,1)$ (resp., $s > 1$), $\g_{\rm tr}$ has a unique extension
to a continuous operator $H^{s+1/2}(\Omega) \rightarrow H^{s}(\G)$
(resp., $H^{s+1/2}(\Omega) \rightarrow H^{1}(\G)$),
see \cite{Costabel_88_BIO,BuffaC_03_EFI}.
We will use the notation $C^\infty(\G) = \g_{\rm tr}(C^\infty(\bar\Omega))$.

Using the introduced Sobolev spaces of scalar functions, we define for $s\ge -1$:
\[
  \bH^s(\Omega) = (H^s(\Omega))^3,\qquad
  \bH^s(\G) = (H^s(\G))^3;
\]
\[
  \bH^s(\G^{(i)}) = (H^s(\G^{(i)}))^2,\quad
  \tilde\bH^s(\G^{(i)}) = (\tH^s(\G^{(i)}))^2,\quad
  1 \le i \le \CI.
\]
If $\G$ is an open surface, then in addition to the above we define the space
\[
  \tilde\bH^s(\G) = (\tH^s(\G))^3,\quad |s| \le 1.
\]
The norms and inner products in all these spaces are defined component-wise and
usual conventions $\bH^0(\Omega) = \bL^2(\Omega)$,
$\bH^0(\G) = \tilde\bH^0(\G) = \bL^2(\G)$,
$\bH^0(\G^{(i)}) = \tilde\bH^0(\G^{(i)}) = \bL^2(\G^{(i)})$ hold.

Now let us introduce the Sobolev spaces of tangential vector fields
defined on $\G$ (see \cite{BuffaC_01_TFI,BuffaC_01_TII,BuffaCS_02_THL}).
We start with the space
\[
  \bL^2_t(\G) := \{\bu \in \bL^2(\G);\; \bu\cdot\bnu = 0\ \hbox{on}\ \G\},
\]
which will be identified with the space of two-dimensional, tangential,
square integrable vector fields. The norm and inner product in this space
will be denoted by $\|\cdot\|_{0,\G}$ and
$\<\cdot,\cdot\> = \<\cdot,\cdot\>_{0,\G}$, respectively.
The similarity of this notation with the one for scalar functions
should not lead to any confusion. Then we define:
\beas
     & \bH_{\;-}^{s}(\G) :=
     \{\bu\in \bL^2_t(\G);\; \bu_i \in \bH^{s}(\G^{(i)}),\ \ 
     1\le i\le \CI\},\quad s \ge 0, &
     \cr\cr
     & \|\bu\|_{\bH_{\;-}^{s}(\G)} :=
     \bigg(\sum\limits_{i=1}^{\CI} \|\bu_i\|_{\bH^{s}(\G^{(i)})}^2\bigg)^{\frac 12}. &
\eeas

Let $\g_{\rm tr}$ be the trace operator (now acting on vector fields),
$\g_{\rm tr}(\bu) = \bu|_{\G}$,
$\g_{\rm tr}: \bH^{s+1/2}(\Omega) \rightarrow \bH^{s}(\G)$ for $s \in (0,1)$,
and let $\g_{\rm tr}^{-1}$ be one of its right inverses.
We will use the ``tangential components trace'' mapping
$\pi_\tau: (C^{\infty}(\bar\Omega))^3 \rightarrow \bL^2_t(\G)$ and
the ``tangential trace'' mapping
$\gamma_\tau: (C^{\infty}(\bar\Omega))^3 \rightarrow \bL^2_t(\G)$,
which are defined as $\bu \mapsto \bnu \times (\bu \times \bnu)|_{\G}$
and $\bu \times \bnu|_{\G}$, respectively.
We will also use the notation $\pi_\tau$ (resp., $\g_\tau$) for the composite
operator $\pi_\tau \circ \g_{\rm tr}^{-1}$ (resp., $\g_\tau \circ \g_{\rm tr}^{-1}$),
which acts on traces. Then we define the spaces
\[
  \bH^{1/2}_{\|}(\G) := \pi_\tau(\bH^{1/2}(\G)),\qquad
  \bH^{1/2}_{\perp}(\G) := \g_\tau(\bH^{1/2}(\G)),
\]
endowed with their operator norms
\beas
     \|\bu\|_{\bH^{1/2}_{\|}(\G)}
     & := &
     \inf_{\bphi \in \bH^{1/2}(\G)} \{\|\bphi\|_{\bH^{1/2}(\G)};\;
     \pi_\tau(\bphi) = \bu\},
     \cr\cr
     \|\bu\|_{\bH^{1/2}_{\perp}(\G)}
     & := &
     \inf_{\bphi \in \bH^{1/2}(\G)} \{\|\bphi\|_{\bH^{1/2}(\G)};\;
     \g_\tau(\bphi) = \bu\}.
\eeas
It has been shown in \cite{BuffaC_01_TFI} that the space $\bH^{1/2}_{\|}(\G)$
(resp., $\bH^{1/2}_{\perp}(\G)$) can be characterised as the space of tangential
vector fields belonging to $\bH_{\;-}^{1/2}(\G)$ and satisfying an appropriate
``weak continuity'' condition for the tangential (resp., normal) component
across each edge $e_{ij}$ of $\G$.

For $s > \frac 12$ we set
\beas
     \bH_{\|}^{s}(\G)
     & := &
     \{\bu \in \bH_{\;-}^{s}(\G);\;
     \bu_i \cdot \btau_{ij} = \bu_j \cdot \btau_{ij}\quad
     \hbox{at each}\ e_{ij}\},
     \cr\cr
     \bH_{\perp}^{s}(\G)
     & := &
     \{\bu \in \bH_{\;-}^{s}(\G);\;
     \bu_i \cdot \btau_{i}^{(j)} = \bu_j \cdot \btau_{j}^{(i)}\quad
     \hbox{at each}\ e_{ij}\}.
\eeas
For any $s > \frac 12$ the spaces $\bH_{\|}^{s}(\G)$ and $\bH_{\perp}^{s}(\G)$
are closed subspaces of $\bH_{\;-}^{s}(\G)$. Finally, for $s \in [0,\frac 12)$
we set
\[
  \bH_{\|}^{s}(\G) = \bH_{\perp}^{s}(\G) := \bH_{\;-}^{s}(\G).
\]

If $\G$ is an open surface, then we also need to define subspaces
of $\bH_{\|}^{s}(\G)$ and $\bH_{\perp}^{s}(\G)$ incorporating boundary
conditions on $\partial\G$ (for tangential and normal components, respectively).
In this case, for a given function $\bu$ on $\G$, we will denote by $\tilde\bu$
the extension of $\bu$ by zero onto a closed Lipschitz polyhedral surface
$\tilde\G \supset \G$. Then we define the spaces
\beas
     \tilde\bH_{\|}^{s}(\G)
     & := &
     \{\bu \in \bH_{\|}^{s}(\G);\;
     \tilde\bu \in \bH_{\|}^{s}(\tilde\G)\},\quad s \ge 0,
     \cr\cr
     \tilde\bH_{\perp}^{s}(\G)
     & := &
     \{\bu \in \bH_{\perp}^{s}(\G);\;
     \tilde\bu \in \bH_{\perp}^{s}(\tilde\G)\},\quad s \ge 0,
\eeas
which are furnished with the norms
\[
  \|\bu\|_{\tilde\bH^{s}_{\|}(\G)} :=
  \|\tilde\bu\|_{\bH^{s}_{\|}(\tilde\G)},\quad
  \|\bu\|_{\tilde\bH^{s}_{\perp}(\G)} :=
  \|\tilde\bu\|_{\bH^{s}_{\perp}(\tilde\G)},\quad s \ge 0.
\]
When considering open and closed surfaces at the same time we use the notation
$\tilde\bH_{\|}^{s}(\G)$, $\tilde\bH_{\perp}^{s}(\G)$, etc. also for closed
surfaces by assuming that $\tilde\bH_{\|}^{s}(\G) = \bH_{\|}^{s}(\G)$,
$\tilde\bH_{\perp}^{s}(\G) = \bH_{\perp}^{s}(\G)$, etc. in this case.
This in particular applies to the following definition of dual spaces.
For $s \in [-1,0)$, the spaces 
$\bH_{\|}^{s}(\G)$, $\tilde\bH_{\|}^{s}(\G)$, $\bH_{\perp}^{s}(\G)$,
$\tilde\bH_{\perp}^{s}(\G)$, and $\bH_{\;-}^{s}(\G)$ are defined as the dual
spaces of $\tilde\bH_{\|}^{-s}(\G)$, $\bH_{\|}^{-s}(\G)$,
$\tilde \bH_{\perp}^{-s}(\G)$, $\bH_{\perp}^{-s}(\G)$, and
$\bH_{\;-}^{-s}(\G)$, respectively (with $\bL^2_t(\G)$ as pivot space).
They are equipped with their natural (dual) norms. Moreover, for any
$s \in (-\frac 12,\frac 12)$ there holds (cf. \cite{Grisvard_92_SBV})
\[
  \tilde\bH_{\|}^{s}(\G) = \bH_{\|}^{s}(\G) =
  \tilde\bH_{\perp}^{s}(\G) = \bH_{\perp}^{s}(\G) = \bH_{\;-}^{s}(\G).
\]

Using the above spaces of tangential vector fields, one can define basic
differential operators on $\G$. The tangential gradient,
$\gradg: H^1(\G) \rightarrow \bL^2_t(\G)$, and the tangential vector curl,
$\bcurlg: H^1(\G) \rightarrow \bL^2_t(\G)$, are defined in the usual way
by localisation to each face $\G^{(i)}$
(see \cite{BuffaC_01_TFI,BuffaC_01_TII} for definitions and properties
of these operators on both closed and open surfaces).
To proceed with extensions of these operators and with their adjoints
we need to distinguish between open and closed surfaces.

Let $\G$ be a closed surface. The adjoint operators of $-\gradg$ and $\bcurlg$
are the surface divergence and the surface scalar curl, respectively:
\be \label{div&curl}
    \divg: \bL^2_t(\G) \rightarrow H^{-1}_{*}(\G),\quad
    \curlg: \bL^2_t(\G) \rightarrow H^{-1}_{*}(\G).
\ee
It has been shown in \cite{BuffaC_01_TII} that $\gradg$ and $\bcurlg$
can be extended to
\[
  \gradg: H^{1/2}(\G) \rightarrow \bH^{-1/2}_{\perp}(\G),\quad
  \bcurlg: H^{1/2}(\G) \rightarrow \bH^{-1/2}_{\|}(\G).
\]
Moreover, they have closed ranges in corresponding spaces. Their
adjoint operators
\[
  \divg: \bH^{1/2}_{\perp}(\G) \rightarrow H^{-1/2}_{*}(\G),\quad
  \curlg: \bH^{1/2}_{\|}(\G) \rightarrow H^{-1/2}_{*}(\G)
\]
are linear continuous and surjective.

Finally, the Laplace-Beltrami operator is defined on $\G$ as follows
\be \label{LB-operator}
    \Delta_\G\,u = \divg(\gradg\, u) = -\,\curlg(\bcurlg\, u)\quad
    \forall u \in H^1(\G).
\ee
One has $\Delta_\G: H^1(\G) \rightarrow H^{-1}_{*}(\G)$,
it is linear continuous and invertible.

If $\G$ is an open surface, then instead of (\ref{div&curl}) there holds
\[
  \divg: \bL^2_t(\G) \rightarrow \tilde H^{-1}(\G),\quad
  \curlg: \bL^2_t(\G) \rightarrow \tilde H^{-1}(\G).
\]
The operators $\gradg$ and $\bcurlg$ again can be extended as follows
(cf. \cite{BuffaC_01_TII}):
\[
  \gradg: \tilde H^{1/2}(\G) \rightarrow \tilde\bH^{-1/2}_{\perp}(\G),\quad
  \gradg: H^{1/2}(\G) \rightarrow \bH^{-1/2}_{\perp}(\G)
\]
and
\[
  \bcurlg: \tilde H^{1/2}(\G) \rightarrow \tilde \bH^{-1/2}_{\|}(\G),\quad
  \bcurlg: H^{1/2}(\G) \rightarrow \bH^{-1/2}_{\|}(\G);
\]
they also have closed ranges in corresponding spaces, and their adjoints
\[
  \divg: \bH^{1/2}_{\perp}(\G) \rightarrow H^{-1/2}(\G),\quad
  \divg: \tilde\bH^{1/2}_{\perp}(\G) \rightarrow \tilde H^{-1/2}(\G)
\]
and
\[
  \curlg: \bH^{1/2}_{\|}(\G) \rightarrow H^{-1/2}(\G),\quad
  \curlg: \tilde\bH^{1/2}_{\|}(\G) \rightarrow \tilde H^{-1/2}(\G)
\]
are linear continuous and surjective.
The Laplace-Beltrami operator $\Delta_\G: H^1(\G) \rightarrow \tH^{-1}(\G)$
is defined as in (\ref{LB-operator}).

We will need the following spaces involving $\Delta_\G$:
\beas
     \CH(\G)
     & := &
     \{u \in H^1(\G)/{\field{C}};\; \Delta_\G\,u \in H^{-1/2}_{*}(\G)\}\quad
     \hbox{if $\G$ is closed,}
     \\[5pt]
     \tilde\CH(\G)
     & := &
     \{u \in H^1(\G)/{\field{C}};\; \Delta_\G\,u \in \tH^{-1/2}(\G)\ \hbox{and}
     \\[3pt]
     &   &
     \quad
     \<\gradg\,u,\gradg\,v\> + \<\Delta_\G\,u,v\> = 0\quad
     \hbox{for all}\ v \in H^1(\G)\}\quad
     \hbox{if $\G$ is an open surface}.
\eeas

Now we can introduce the spaces which appear when dealing with the EFIE on $\G$.
First, we set
\[
  \bH^s(\divg,\G) :=
  \{\bu \in \bH^s_{\|}(\G);\; \divg\,\bu \in H^s(\G)\},\quad s \in [-1/2,0]
\]
(here, $\G$ is either a closed or an open surface).
If $\G$ is an open surface, then we will also use the space
\[
  \tilde\bH^s(\divg,\G) :=
  \{\bu \in \tilde\bH^{s}_{\|}(\G);\;
  \divg\,\bu \in \tH^{s}(\G)\},\quad s \in [-1/2,0].
\]
The spaces $\bH^s(\divg,\G)$ and $\tilde\bH^s(\divg,\G)$ are equipped with their graph
norms $\|\cdot\|_{\bH^s(\divg,\G)}$ and $\|\cdot\|_{\tilde\bH^s(\divg,\G)}$,
respectively. For $s=0$ we drop the superscript and for open surfaces also
the tilde in the above notation,
$\bH^0(\divg,\G) = \tilde\bH^0(\divg,\G) = \bH(\divg,\G)$.

On open surfaces, one also needs the spaces incorporating homogeneous
boundary conditions for the trace of the normal component on $\partial\G$.
By $\bH_0(\divg,\G)$ (resp., $\tilde\bH^{-1/2}_0(\divg,\G)$) we denote the
subspace of elements $\bu \in \bH(\divg,\G)$ (resp., $\bu \in \tilde\bH^{-1/2}(\divg,\G)$)
such that for all $v \in C^{\infty}(\G)$ there holds
\be \label{X_property}
    \<\bu,\gradg v\> + \<\divg\,\bu,v\> = 0.
\ee
We note that if $\bu \in \tilde\bH^{-1/2}_0(\divg,\G)$,
then identity (\ref{X_property}) holds for any $v\in H^{3/2}(\G)$ by density.
In particular, $\tilde\bH^{-1/2}_0(\divg,\G)$ is a closed subspace of
$\tilde\bH^{-1/2}(\divg,\G)$. To join the notation for open and closed surfaces,
we will write
\[
  \bX^0 = \bH(\divg,\G),\quad \bX = \bH^{-1/2}(\divg,\G)\quad
  \hbox{if $\G$ is a closed surface},
\]
\[
  \bX^0 = \bH_0(\divg,\G),\quad \bX = \tilde\bH^{-1/2}_0(\divg,\G)\quad
  \hbox{if $\G$ is an open surface}.
\]
In both cases the norm in the space $\bX$ will be denoted as $\|\cdot\|_{\bX}$.

\subsection{Some equivalent norms and corresponding inner products} \label{sec_norms}

In this subsection we consider the Sobolev spaces $H^s$ and $\tH^s$
on the reference square $K$ for $s = \pm\frac 12$.
We will derive expressions for norms which are
equivalent to those introduced in Section~\ref{sec_spaces}.
We note that all results of this subsection are valid also for $K$
being the equilateral reference triangle.
First, let us introduce some notation.

\begin{enumerate}

\item[$1^{\circ}$.]
We denote by $\Omega$ the cube $\Omega = K \times (0,1)$.
Thus $\partial\Omega = \cup_{i=1}^{6} \bar\G_i$.
Let $K = \G_1 = \{(x_1,x_2,0);\; (x_1,x_2) \in K\}$,
$\G_{6} = \{(x_1,x_2,1);\; (x_1,x_2) \in K\}$, and denote
$\tilde K = \partial\Omega \backslash \bar\G_{6}$.
Note that $\tilde K$ is an open surface.
We will use the standard notation for the gradient $\nabla$ and for the Laplace
operator $\Delta$, both acting on scalar functions of three variables.

\item[$2^{\circ}$.]
Given $u \in H^{-1/2}(K)$, we denote by $\tilde u_K$ the solution of the mixed
problem: find $\tilde u_K \in H^1(\Omega)$ such that
\[
  \Delta \tilde u_K = 0\ \hbox{in $\Omega$},\quad
  \hbox{$\frac{\partial\tilde u_K}{\partial n} = u$ on $K$},\quad
  \tilde u_K = 0\ \hbox{on $\partial\Omega \backslash K$}.
\]
If $u \in H^{-1/2}(\tilde K)$, then we will use the same notation as above
with $K$ replaced by $\tilde K$.

\item[$3^{\circ}$.]
Given $u \in H^{1/2}(\partial\Omega)$, we denote by $\tilde{\tilde u}$
its harmonic extension, i.e., the solution of the Dirichlet problem:
find $\tilde{\tilde u} \in H^1(\Omega)$ such that
\be \label{Dir_ext}
    \Delta \tilde{\tilde u} = 0\ \hbox{in $\Omega$},\quad
    \tilde{\tilde u} = u\ \hbox{on $\partial\Omega$}.
\ee

\item[$4^{\circ}$.]
Given $u \in \tH^{1/2}(K)$, we denote by $u^{\circ}$ the extension
of $u$ by zero onto $\partial\Omega$.
Thus, $u^{\circ} \in H^{1/2}(\partial\Omega)$.

\end{enumerate}

We make use of standard definitions for the norm and the semi-norm in $H^1(\Omega)$:
\[
   \|u\|_{H^1(\Omega)} =
   \left(\|u\|_{0,\Omega}^2 + |u|_{H^1(\Omega)}^2\right)^{1/2},\quad
   |u|_{H^1(\Omega)} = \|\grad u\|_{0,\Omega}.
\]
Since $H^{1/2}(\partial\Omega)$ is the trace space of $H^1(\Omega)$,
the norm and the semi-norm in $H^{1/2}(\partial\Omega)$ can be equivalently
written as follows
\bea
    \|u\|_{H^{1/2}(\partial\Omega)}
    & \simeq &
    \stack{U|_{\partial\Omega}=u}{\inf_{U \in H^1(\Omega)}} \|U\|_{H^1(\Omega)},
    \nonumber
    \\
    |u|_{H^{1/2}(\partial\Omega)}
    & \simeq &
    \stack{U|_{\partial\Omega}=u}{\inf_{U \in H^1(\Omega)}} |U|_{H^1(\Omega)} =
    \|\grad \tilde{\tilde u}\|_{0,\Omega}.
    \label{1/2-seminorm}
\eea
Now we can define equivalent norms in $\tilde H^{1/2}(K)$ and $H^{1/2}(K)$:
\bea
    \|u\|_{\tH^{1/2}(K)}
    & \simeq &
    |u^\circ|_{H^{1/2}(\partial\Omega)} \simeq
    \Big\|\grad \widetilde{\widetilde{u^\circ}}\Big\|_{0,\Omega},
    \label{tilde_1/2-norm}
    \\
    \|u\|_{H^{1/2}(K)}
    & \simeq &
    \stack{U|_{K}=u}{\inf_{U \in \tilde H^{1/2}(\tilde K)}} \|U\|_{\tH^{1/2}(\tilde K)},
    \label{1/2-norm}
\eea
where $\|\cdot\|_{\tH^{1/2}(\tilde K)}$ is defined as in (\ref{tilde_1/2-norm}),
because $\tilde K$ is an open surface.

From (\ref{tilde_1/2-norm}) one can easily derive the expression for the
corresponding $\tH^{1/2}(K)$-inner product.
In fact, applying the parallelogram law twice, integrating by parts, and
recalling notations $3^\circ$, $4^\circ$, we find
(see also \cite{DemkowiczB_03_pIE})
\bea
    \<u,v\>_{\tH^{1/2}(K)}
    & = &
    \Big\<\grad \ttu,\grad \ttv\Big\>_{0,\Omega} =
    \Big\<\frac{\partial\ttu}{\partial n},\ttv\Big\>_{0,\partial\Omega} =
    \nonumber
    \\[5pt]
    & = &
    \Big\<\frac{\partial\ttu}{\partial n},v\Big\>_{0,K} =
    \Big\<u,\frac{\partial\ttv}{\partial n}\Big\>_{0,K}\quad
    \forall u,v \in \tH^{1/2}(K).
    \label{tilde_1/2-ip}
\eea

The space $H^{-1/2}(K)$ is the dual space of $\tH^{1/2}(K)$.
We prove the following result regarding an equivalent norm in $H^{-1/2}(K)$.

\begin{lemma} \label{lm_-1/2-norm}
For any $u \in H^{-1/2}(K)$ there holds
\be \label{-1/2-norm}
    \|u\|_{H^{-1/2}(K)} \simeq \|\grad \tilde u_K\|_{0,\Omega}.
\ee
The $H^{-1/2}$-inner product corresponding to the norm on the
right-hand side of {\rm (\ref{-1/2-norm})} reads as
\be \label{-1/2_ip}
    \<u,v\>_{H^{-1/2}(K)} = \<u,\tilde v_K\>_{0,K} = \<\tilde u_K,v\>_{0,K}\quad
    \forall u,v \in H^{-1/2}(K).
\ee
\end{lemma}

\begin{proof}
Using notations $2^\circ-4^\circ$,
we integrate by parts to obtain for any $u \in H^{-1/2}(K)$ and
any $v \in \tH^{1/2}(K)$
\[
      \Big\<\grad \tilde u_K,\grad \ttv\Big\>_{0,\Omega} =
      \Big\<\frac{\partial\tilde u_K}{\partial n},\ttv\Big\>_{0,\partial\Omega} =
      \Big\<\frac{\partial\tilde u_K}{\partial n},\ttv\Big\>_{0,K} +
      \Big\<\frac{\partial\tilde u_K}{\partial n},\ttv\Big\>_{0,\partial\Omega\backslash K} =
      \<u,v\>_{0,K}.
\]
Hence, we find from (\ref{dual_norm}) and (\ref{tilde_1/2-norm})
\be \label{-1/2-norm_aux1}
    \|u\|_{H^{-1/2}(K)} =
    \sup_{0\not=v\in\tH^{1/2}(K)}
    {\Big|\Big\<\grad \tilde u_K,\grad \ttv\Big\>_{0,\Omega}\Big|
    \over{\|v\|_{\tH^{1/2}(K)}}} \simeq
    \sup_{0\not=v\in\tH^{1/2}(K)}
    {\Big|\Big\<\grad \tilde u_K,\grad \ttv\Big\>_{0,\Omega}\Big|
    \over{\Big\|\grad \ttv\Big\|_{0,\Omega}}}.
\ee
Let $w := \tilde u_K|_K$. One has $w \in \tH^{1/2}(K)$
because $\tilde u_K = 0$ on $\partial\Omega\backslash K$. Moreover,
$w^\circ = \tilde u_K|_{\partial\Omega}$ and, due to the uniqueness of the
solution to the Dirichlet problem (\ref{Dir_ext}), we conclude that
$\widetilde{\widetilde{w^\circ}} = \tilde u_K$. Therefore,
\be \label{-1/2-norm_aux2}
    \sup_{0\not=v\in\tH^{1/2}(K)}
    {\Big|\Big\<\grad \tilde u_K,\grad \ttv\Big\>_{0,\Omega}\Big|
    \over{\Big\|\grad \ttv\Big\|_{0,\Omega}}} \ge
    {\Big|\Big\<\grad \tilde u_K,\grad \widetilde{\widetilde{w^\circ}}\Big\>_{0,\Omega}\Big|
    \over{\Big\|\grad \widetilde{\widetilde{w^\circ}}\Big\|_{0,\Omega}}} =
    \|\grad \tilde u_K\|_{0,\Omega}.
\ee
On the other hand, it is easy to see that
\be \label{-1/2-norm_aux3}
    \sup_{0\not=v\in\tH^{1/2}(K)}
    {\Big|\Big\<\grad \tilde u_K,\grad \ttv\Big\>_{0,\Omega}\Big|
    \over{\Big\|\grad \ttv\Big\|_{0,\Omega}}} \le
    \|\grad \tilde u_K\|_{0,\Omega}.
\ee
Now (\ref{-1/2-norm}) immediately follows from 
(\ref{-1/2-norm_aux1})--(\ref{-1/2-norm_aux3}).

Using (\ref{-1/2-norm}) together with the parallelogram law we find
\[
  \<u,v\>_{H^{-1/2}(K)} =
  \Big\<\grad \tilde u_K,\grad \tilde v_K\Big\>_{0,\Omega}\quad
  \forall u,v \in H^{-1/2}(K).
\]
Hence, integrating by parts and using notation~$2^\circ$, we derive (\ref{-1/2_ip}).
\end{proof}

The following lemma states an analogous result for the space $\tH^{-1/2}(K)$
which is the dual space of $H^{1/2}(K)$.

\begin{lemma}
For any $u \in \tH^{-1/2}(K)$ there holds
\be \label{tilde-1/2-norm}
    \|u\|_{\tH^{-1/2}(K)} \simeq
    \Big\|\grad \widetilde{(u^\circ)}_{\tilde K}\Big\|_{0,\Omega}.
\ee
The $\tH^{-1/2}$-inner product corresponding to the norm on the
right-hand side of {\rm (\ref{tilde-1/2-norm})} reads as
\be \label{tilde_-1/2-ip}
    \<u,v\>_{\tH^{-1/2}(K)} =
    \Big\<u,\widetilde{(v^\circ)}_{\tilde K}\Big\>_{0,K} =
    \Big\<\widetilde{(u^\circ)}_{\tilde K},v\Big\>_{0,K}\quad
    \forall u,v \in \tH^{-1/2}(K).
\ee
\end{lemma}

\begin{proof}
Let $u \in \tH^{-1/2}(K)$.
Then $u^\circ \in \tH^{-1/2}(\tilde K) \subset H^{-1/2}(\tilde K)$.
Using (\ref{dual_norm}) and (\ref{1/2-norm}) we have
\beas
     \|u^\circ\|_{H^{-1/2}(\tilde K)}
     & = &
     \sup_{0\not=w\in\tH^{1/2}(\tilde K)}
     {|\<u^\circ,w\>_{0,\tilde K}|\over{\|w\|_{\tH^{1/2}(\tilde K)}}} =
     \sup_{0\not=w\in\tH^{1/2}(\tilde K)}
     {|\<u,w\>_{0,K}|\over{\|w\|_{\tH^{1/2}(\tilde K)}}}
     \\[5pt]
     & = &
     \sup_{0\not=v\in H^{1/2}(K)}\,\,
     \stack{V|_{K}=v}{\sup_{V \in \tH^{1/2}(\tilde K)}}
     {|\<u,V\>_{0,K}|\over{\|V\|_{\tH^{1/2}(\tilde K)}}} =
     \sup_{0\not=v\in H^{1/2}(K)}
     \frac{|\<u,v\>_{0,K}|}
     {\displaystyle{
                    \stack{V|_{K}=v}{\inf_{V\in\tH^{1/2}(\tilde K)}} \|V\|_{\tH^{1/2}(\tilde K)}
                   }
     }
     \\[3pt]
     & \simeq &
     \sup_{0\not=v\in H^{1/2}(K)}
     {|\<u,v\>_{0,K}|\over{\|v\|_{H^{1/2}(K)}}} =
     \|u\|_{\tH^{-1/2}(K)}.
\eeas
Hence, using (\ref{-1/2-norm}) with $u$ replaced by $u^\circ$ and with $K$
replaced by $\tilde K$, we prove (\ref{tilde-1/2-norm}):
\[
  \|u\|_{\tH^{-1/2}(K)} \simeq \|u^\circ\|_{H^{-1/2}(\tilde K)} \simeq
  \Big\|\grad \widetilde{(u^\circ)}_{\tilde K}\Big\|_{0,\Omega}\quad
  \forall u \in \tH^{-1/2}(K).
\]
Then, applying the parallelogram law, integrating by parts, and
making use of notations $2^\circ,\;4^\circ$, we derive (\ref{tilde_-1/2-ip}).
\end{proof}

\begin{remark} \label{rem_edge-norms}
The same arguments as above can be used to find equivalent norms and corresponding
inner products in the Sobolev spaces on edges $\ell_i \subset \partial K$.
In particular, we will need the $\tH^{1/2}(\ell_i)$-norm and corresponding
inner product. Using the notation analogous to $3^\circ$ and $4^\circ$, we have
(cf. {\rm (\ref{tilde_1/2-norm}), (\ref{tilde_1/2-ip})})
\beas
     \|u\|_{\tH^{1/2}(\ell_i)}
     & \simeq &
     \Big\|\grad \widetilde{\widetilde{u^\circ}}\Big\|_{0,K}\quad
     \forall u \in \tH^{1/2}(\ell_i),
     \\
     \<u,v\>_{\tH^{1/2}(\ell_i)}
     & = &
     \Big\<\frac{\partial\ttu}{\partial n},v\Big\>_{0,\ell_i} =
     \Big\<u,\frac{\partial\ttv}{\partial n}\Big\>_{0,\ell_i}\quad
     \forall u,v \in \tH^{1/2}(\ell_i).
\eeas
\end{remark}

The next lemma states the fact that for a constant function $v$ in
(\ref{tilde_-1/2-ip}) the $\tH^{-1/2}(K)$-inner product reduces to the
$L^2(K)$-inner product.

\begin{lemma} \label{lm_tilde_-1/2_ip}
For any $u \in \tH^{-1/2}(K)$ there holds
\[
  \<u,1\>_{\tH^{-1/2}(K)} = \<u,1\>_{0,K}.
\]
\end{lemma}

\begin{proof}
We have by (\ref{tilde_-1/2-ip})
\be \label{lm_tilde1}
    \<u,1\>_{\tH^{-1/2}(K)} = \<u,\varphi|_K\>_{0,K},
\ee
where $\varphi(x)$ ($x = (x_1,x_2,x_3) \in \Omega = K \times (0,1)$)
solves the following mixed problem (see (\ref{tilde_-1/2-ip})
and notations ${1^\circ}$, ${2^\circ}$, $4^\circ$):
find $\varphi \in H^1(\Omega)$ such that
\[
  \Delta \varphi = 0\ \hbox{in $\Omega$},\quad
  \hbox{$\frac{\partial\varphi}{\partial n} = 1$ on $\G_1 = K$},\quad
  \hbox{$\frac{\partial\varphi}{\partial n} = 0$ on $\G_i$ ($i = 2,\ldots,5$)},\quad
  \varphi = 0\ \hbox{on $\G_{6}$}.
\]
It is easy to see that $\varphi = 1 - x_3$.
Then $\varphi|_K = \varphi|_{x_3=0} = 1$ and the assertion follows
from (\ref{lm_tilde1}).
\end{proof}

\subsection{Auxiliary lemmas} \label{sec_aux}

The Laplace-Beltrami operator $\Delta_\G$ will be a useful tool in our analysis.
The following two lemmas establish its regularity separately on closed and open
piecewise plane Lipschitz surfaces. For proofs we refer to
\cite[Theorem~8]{BuffaCS_02_BEM} and \cite[Proposition~4.11]{BuffaC_03_EFI},
respectively.

\begin{lemma} \label{lm_LB_reg1}
Let $\G$ be a closed Lipschitz polyhedral surface, $\psi \in H^{s}_{*}(\G)$
for $s \ge -1$, and let $\phi \in H^1(\G)/{\field{C}}$ be the unique
solution to the problem
\[
  \<\gradg\,\phi, \gradg\,\tilde\phi\> = \<\psi,\tilde\phi\>\quad
  \forall \tilde\phi \in H^1(\G)/{\field{C}}.
\]
Then $\phi \in H^{1+r}(\G)$ and
\[
  \|\phi\|_{H^{1+r}(\G)/{\field{C}}} \le C\, \|\psi\|_{H^s(\G)}
\]
for any $r < \min\,\{s^*,s+1\}$, where $s^* > 0$ depends on the geometry
of $\G$ in neighbourhoods of its vertices.
\end{lemma}

\begin{lemma} \label{lm_LB_reg2}
Let $\G$ be an open piecewise plane Lipschitz surface.
Let $\psi \in H^{s}_{*}(\G)$, $s > -\frac 12$
(resp.,  $\psi \in \tH^{s}(\G)$, $-1 \le s \le -\frac 12$,
$\<\psi,1\> = 0$) and
$\phi \in H^1(\G)/{\field{C}}$ be the unique solution to the problem
\[
  \<\gradg\,\phi, \gradg\,\tilde\phi\> = \<\psi,\tilde\phi\>\quad
  \forall \tilde\phi \in H^1(\G)/{\field{C}}.
\]
Then $\phi \in H^{1+r}(\G)$ and
\[
  \|\phi\|_{H^{1+r}(\G)/{\field{C}}} \le C\, \|\psi\|_{H^s(\G)}\quad
  \hbox{(resp.,}\quad
  \|\phi\|_{H^{1+r}(\G)/{\field{C}}} \le C\, \|\psi\|_{\tH^s(\G)})
\]
for any $r < \min\,\{s^*,s+1\}$, where $s^* > 0$ depends on the geometry
of $\G$ in neighbourhoods of all vertices of $\bar\G$.
\end{lemma}

\begin{remark} \label{rm_LB_reg}
If $s+1 < s^*$ in the above two lemmas for given $\G$ and $s$,
then both results are valid for $0 \le r \le s+1$.
\end{remark}

In the following lemma we formulate some useful properties of the Piola transform.

\begin{lemma} \label{lm_scaling}
Let $K^h$ and $K$ be two open subsets of ${\field{R}}^2$ such that $K^h = T(K)$,
where $T$ is a sufficiently smooth one-to-one mapping with a sufficiently
smooth inverse $T^{-1}:\,K^h \rightarrow K$.
Assume that $\diam\, K^h\simeq \rho_{K^h}\simeq h$, $\diam\, K\simeq \rho_{K}\simeq 1$,
and the mappings $T,\,T^{-1}$ satisfy the same relations as in
{\rm (\ref{Jacobian})}--{\rm (\ref{mesh-reg})}.
Let $\hat\varphi$ and $\hat\bq$ be a scalar function and a vector field, respectively,
defined on $K$, and let $\varphi = \hat\varphi \circ T^{-1}$, $\bq = \CM(\hat\bq)$
be defined on $K^h$ (here, $\CM$ is the Piola transform associated with $T$, see
{\rm (\ref{Piola})}). Then
\bea
    \<\varphi,\div\,\bq\>_{0,K^h}
    & = &
    \<\hat\varphi,\div\,\hat\bq\>_{0,K},
    \label{Piola_1}
    \\[5pt]
    \|\bq\|_{0,K^h}
    & \simeq &
    \|\hat\bq\|_{0,K}
    \label{Piola_2}
\eea
if $\hat\varphi \in L^2(K)$ and $\hat\bq \in \bH(\div,K)$.
Moreover, for $s \in [0,1]$, there holds
\be \label{Piola_3}
    \|\hat\bq\|_{\bH^s(K)} \le C\, \|\bq\|_{\bH^s(K^h)}
\ee
if $\bq \in \bH^s(K^h)$;
\be \label{Piola_4}
    \|\div\,\hat\bq\|_{\tH^{-s}(K)} \le
    C\,h^{1-s}\, \|\div\,\bq\|_{\tH^{-s}(K^h)}
\ee
if $\div\,\bq \in \tH^{-s}(K^h)$;
\be \label{Piola_5}
    \|\div\,\hat\bq\|_{H^{-s}(K)} \simeq
    h^{1-s}\, \|\div\,\bq\|_{H^{-s}(K^h)}
\ee
if $\div\,\bq \in H^{-s}(K^h)$.
\end{lemma}

\begin{proof}
Relations (\ref{Piola_1}), (\ref{Piola_2}) are well-known
(see, e.g., Lemmas~1.5,~1.6 in Chapter~III of~\cite{BrezziF_91_MHF}).
If $\bq \in \bH^1(K^h)$ then (cf. \cite{Thomas_thesis})
\be \label{Piola_6}
    \|\hat\bq\|_{\bH^1(K)} \le C\, \|\bq\|_{\bH^1(K^h)}.
\ee
Then (\ref{Piola_3}) follows from (\ref{Piola_2}) and (\ref{Piola_6}) by interpolation.

In order to prove (\ref{Piola_4}) and (\ref{Piola_5}) we use (\ref{Piola_1})
and the standard scaling argument for scalar functions. For instance, in the former
case we have
\beas
     \|\div\,\hat\bq\|_{\tH^{-s}(K)}
     & = &
     \sup_{0 \not= \hat\varphi \in H^{s}(K)}
     {\<\div\,\hat\bq,\hat\varphi\>_{0,K} \over{\|\hat\varphi\|_{H^{s}(K)}}}
     \\[5pt]
     & \le &
     C\,\sup_{0 \not= \varphi \in H^{s}(K^h)}
     {\<\div\,\bq,\varphi\>_{0,K^h} \over{h^{-(1-s)}\,\|\varphi\|_{H^{s}(K^h)}}}
     = C\, h^{1-s}\, \|\div\,\bq\|_{\tH^{-s}(K^h)}.
\eeas
The proof of (\ref{Piola_5}) is analogous.
\end{proof}

The following lemma states the inverse inequality for polynomials on $K$.
We refer to \cite{Heuer_01_ApS} for a proof.

\begin{lemma} \label{lm_inverse}
Let $v_p \in \CP_p(K)$. Then for any $s,\,r \in [-1,1]$ with
$s \le r$ there holds
\[
  \|v\|_{H^r(K)} \le C\, p^{2(r-s)}\, \|v\|_{H^s(K)},
\]
where $C$ is a positive constant independent of $p$.
\end{lemma}

\section{Decompositions} \label{sec_dec}
\setcounter{equation}{0}

The main tool in the analysis of the EFIE is the Helmholtz decomposition
of the energy space $\bX$. It is used to prove an inf-sup condition for
the electric field integral operator and to establish
the unique solvability of the EFIE on $\G$
(see, e.g., \cite{BuffaCS_02_BEM,BuffaC_03_EFI}).
The following statement establishes the Helmholtz decomposition of $\bX$
on a (closed or open) Lipschitz polyhedral surface $\G$.
This result has been proved in \cite[Theorems~5.1 and~6.4]{BuffaC_01_TII}
(for open surfaces see also \cite[Section~2.4]{BuffaC_03_EFI}).

\begin{theorem} \label{thm_H-dec}
Let
\beas
     \bW
     & = &
     \{\bw \in \bX;\; \divg\,\bw = 0\},
     \\
     \bV
     & = &
     \{\bv \in \bX;\; \<\bv,\bw\> = 0\ \ \forall \bw \in \bW \cap \bL^2_t(\G)\}.
\eeas
Then there holds
\be \label{H-dec}
    \bX = \bV \oplus \bW.
\ee
Furthermore, $\bV$ and $\bW$ are closed subspaces of $\bX$, and they can be
written as
\[
  \bV = \gradg \CH(\G),\qquad \bW = \bcurlg\,(H^{1/2}(\G)/\field{C})
\]
if $\G$ is a closed surface, and
\[
  \bV = \gradg \tilde\CH(\G),\qquad \bW = \bcurlg\, \tH^{1/2}(\G)
\]
if $\G$ is an open surface.
\end{theorem}

In this paper we discretise the EFIE by the $hp$-version of the Galerkin BEM
based on the sequence of the RT-subspaces $\bX_{hp} \subset \bX$
(see (\ref{Xp}), (\ref{BEM})). 
To prove the well-posedness of (\ref{BEM}) (see Theorem~\ref{thm_solve})
we follow \cite{BuffaC_03_EFI,BuffaHvPS_03_BEM} and consider
$\bL^2_t(\G)$-orthogonal discrete decompositions of $\bX_{hp}$ mimicking
the Helmholtz decomposition of $\bX$:
\be \label{d-dec}
    \bX_{hp} = \bV_{hp} \oplus \bW_{hp},
\ee
where
\bea
    \bW_{hp}
    & := &
    \{\bw \in \bX_{hp};\; \divg\,\bw = 0\},
    \nonumber
    \\
    \bV_{hp}
    & := &
    \{\bv \in \bX_{hp};\; \<\bv,\bw\> = 0\ \ 
    \forall \bw \in \bW_{hp}\}.
    \label{V_p}
\eea
It is easy to see that $\bW_{hp} \subset \bW$, however, in general,
$\bV_{hp} \not\subset \bV$. That is why the discrete inf-sup condition
(and thus, the unique solvability of (\ref{BEM}) and quasi-optimal convergence of
the BEM) cannot be deduced by standard arguments, which are usually applied to
conforming Galerkin discretisations of coercive variational problems.

Sufficient conditions to establish the well-posedness of the Galerkin BEM
applied to problem (\ref{bie_var}) were found in \cite{BuffaC_03_EFI}.
It turns out that it is enough to prove that discrete decompositions (\ref{d-dec})
are in some sense close to the Helmholtz decomposition (\ref{H-dec})
of $\bX$ when the dimension of the discrete space tends to infinity.
The abstract formulation of this approach is given in the next theorem
(here, we quote \cite[Theorem~4.1]{BuffaHvPS_03_BEM}, see also
Proposition~4.1, Corollary~4.2, and Theorem~4.5 in \cite{BuffaC_03_EFI}).

\begin{theorem} \label{thm_abstract}
Let $\{\bX_n\}_n$ be a sequence of closed subspaces $\bX_n \subset \bX$ with
decompositions $\bX_n = \bV_n \oplus \bW_n$ which are stable with respect to complex
conjugation and which satisfy the following assumptions:
\begin{enumerate}
\item[{\rm (A1)}]
the family $\{\bX_n\}_n$ is dense in the space $\bX$, namely
\[
  \overline{\bigcup_n \bX_n} = \bX;
\]
\item[{\rm (A2)}]
the spaces $\bV_n$ and $\bW_n$ are such that $\bW_n \subset \bW$ and
\be \label{A2}
    \sup_{\bv_n \in \bV_n\setminus\{\bzero\}}\,\inf_{\bv \in \bV}
    \frac{\|\bv_n - \bv\|_{\bX}}{\|\bv_n\|_{\bX}}\,\rightarrow\,0\ \ 
    \hbox{as}\ \ n\rightarrow\infty.
\ee
\end{enumerate}
Then there exists $n_0$ such that for all $\boldf \in \bX'$ and $n \ge n_0$ the Galerkin
system
\[
  a(\bu_n,\bv) = \<\boldf, \bv\> \quad \forall \bv \in \bX_{n}
\]
has a unique solution $\bu_n \in \bX_n$ which converges quasi-optimally, i.e.,
\[
  \|\bu - \bu_{n}\|_{\bX} \le
  C \inf\{\|\bu - \bv\|_{\bX};\; \bv\in \bX_{n}\},
\]
where $\bu \in \bX$ is the solution of {\rm (\ref{bie_var})}.
\end{theorem}

It has been proved in \cite{BuffaC_03_EFI, BuffaHvPS_03_BEM} that
$\bL^2_t(\G)$-orthogonal discrete decompositions mimicking the Helmholtz decomposition
of the space $\bX$ satisfy assumptions (A1), (A2) of Theorem~\ref{thm_abstract} with
respect to the mesh parameter $h$, i.e., in the framework of the
$h$-version of the BEM for the EFIE. To prove this result for the $hp$-version
on quasi-uniform meshes, we will need an $\tilde\bH^{-1/2}(\div)$-conforming
$p$-interpolation operator, which is introduced and analysed in the next~section.

\section{Interpolation operators} \label{sec_inter}
\setcounter{equation}{0}

The main purpose of this section is to introduce and analyse a new
$\tilde\bH^{-1/2}(\div)$-conforming $p$-interpolation operator.
This operator is necessary to deal with
low regular vector fields, such as gradients of solutions to boundary value
problems for the Laplace-Beltrami operator on polyhedral surfaces (see the regularity
results of Lemmas~\ref{lm_LB_reg1} and~\ref{lm_LB_reg2}). We will construct
the $\tilde\bH^{-1/2}(\div)$-conforming interpolation operator by employing
the $\tH^{-1/2}$-projection for the divergence term. For this operator we then
prove quasi-stability with respect to polynomial degrees and commutativity
with the $\tH^{-1/2}$-projector.

In this section we use standard differential operators $\div$, $\curl$ and
$\grad$, $\bcurl$ acting on 2D vector fields and scalar functions, respectively.
First, let us recall some known interpolation operators acting
on scalar functions and vector fields on $K$.
Let $\Pi^0_p:\, L^2(K) \rightarrow \CP_p(K)$ be the standard
$L^2$-projection onto the set of polynomials $\CP_p(K)$.
We will also use the $\tH^{-1/2}$-projector onto $\CP_p(K)$
denoted by $\Pi^{-1/2}_p:\, \tH^{-1/2}(K) \rightarrow \CP_p(K)$
and satisfying
\[
  \<u - \Pi^{-1/2}_p\,u,v\>_{-\frac 12,K} = 0\quad
  \forall v \in \CP_p(K).
\]
Here and below $\<\cdot,\cdot\>_{-\frac 12,K}$ denotes the $\tH^{-1/2}(K)$-inner product
(see (\ref{tilde_-1/2-ip})).

In \cite{DemkowiczB_03_pIE} two projection-based
interpolation operators have been introduced and analysed. These are
the $H^1$-conforming interpolation operator
$\Pi^{1}_p:\, H^{1+r}(K) \rightarrow \CP_p(K)$ and
the $\bH(\curl)$-conforming interpolation operator
$\Pi^{\curl}_p:\, \bH^r(K) \cap \bH(\curl,K) \rightarrow \bCP^{\rm Ned}_p(K)$
(here, $r > 0$ in both cases and $\bCP^{\rm Ned}_{p}(K) = \CP_{p-1,p}(K) \times \CP_{p,p-1}(K)$
is the $\bH(\curl)$-conforming (first) N{\' e}d{\' e}lec space of degree $p$).
Later, in \cite{BoffiCDD_06_Dhp} these operators were employed
to prove the discrete compactness property for $hp$ adaptive rectangular edge
finite elements. Due to the isomorphism of
the curl and the $\div$ operator in 2D (and, as a consequence, the isomorphism
of the N{\' e}d{\' e}lec elements of the first type and the RT elements),
one can use the results of \cite{DemkowiczB_03_pIE,BoffiCDD_06_Dhp}
related to the operator $\Pi^{\curl}_p$ in the $\bH(\div)$-settings.

We will denote by $\Pi^{\div,0}_p$ the corresponding $\bH(\div)$-conforming
projection-based interpolation operator. Then for $r > 0$ the following diagram commutes
(see Proposition~3 in~\cite{DemkowiczB_03_pIE}):
\be \label{deRham}
    \ba{ccccc}
    H^{1+r}(K)                & \stackrel{\bcurl}{\longrightarrow} &
    \bH^r(K) \cap \bH(\div,K) & \stackrel{\div}{\longrightarrow}   & L^2(K)
    \cr
    \quad\left\downarrow{\Large\strut}\right.\,\Pi^1_p         &
                                                               &
    \qquad\left\downarrow{\Large\strut}\right.\,\Pi^{\div,0}_p &
                                                               &
    \qquad\left\downarrow{\Large\strut}\right.\,\Pi^{0}_{p-1}
    \cr
    \CP_{p}(K)          & \stackrel{\bcurl}{\longrightarrow} &
    \bCP^{\rm RT}_p(K)  & \stackrel{\div}{\longrightarrow}   & \CP_{p-1}(K).
    \ea
\ee
Furthermore, we use the above mentioned isomorphisms to reformulate
the following two results from~\cite{BoffiCDD_06_Dhp}.

\begin{lemma} \label{lm_inf-sup}
{\rm \cite[Theorem~4]{BoffiCDD_06_Dhp}}
Let $\bA_p = \bCP^{\rm RT,0}_p(K)$ and
$\bB_p = \bcurl\,\CP_p^0(K) \oplus \grad\,\div\,\bCP^{\rm RT,0}_p(K)$.
Then the following stability condition holds
\[
  \inf_{\bolda \in \bA_p} \sup_{\boldb \in \bB_p}
  {\<\bolda,\boldb\>_{0,K}\over{\|\bolda\|_{0,K}\,\|\boldb\|_{0,K}}} =
  \left({2(2p+1) \over {(p+1)(p+2)}}\right)^{1/2} =
  O(p^{-1/2}).
\]
\end{lemma}

This lemma and the definition of the interpolation operator $\Pi^{\div,0}_p$
imply the following $\bL^2$-stability result for the $p$-version.

\begin{lemma} \label{lm_L^2-stab}
{\rm \cite[Theorem~7]{BoffiCDD_06_Dhp}}
Let $\bu \in \bH^{r}(K) \cap \bH_0(\div,K)$, $r > 0$, be a curl-free bubble function on $K$,
and let $\bu_p^{\div} := \Pi_p^{\div,0} \bu$. Then $\bu_p^{\div} \in \bCP^{\rm RT,0}_p(K)$
is discrete curl-free (i.e., $\<\bu_p^{\div},\bcurl\,\phi\>_{0,K} = 0$
for any $\phi \in \CP_p^0(K))$ and there holds
\[
  \|\bu - \bu_p^{\div}\|_{0,K} \le
  C\,p^{1/2}\,\inf_{\bq_p \in \bCP^{\rm RT,0}_p(K)} \|\bu - \bq_p\|_{0,K}.
\]
\end{lemma}

In the next lemma we estimate the error of the best $\bL^2$-approximation
of low regular vector bubble functions by RT-elements of degree $p$.

\begin{lemma} \label{lm_L^2-approx}
Let $\bu \in \bH^{r}(K) \cap \bH_0(\div,K)$, $r \in [0,1]$, be a general bubble
function on $K$. Then
\[
  \inf_{\bq_p \in \bCP^{\rm RT,0}_p(K)} \|\bu - \bq_p\|_{0,K} \le
  C\,p^{-r}\,\|\bu\|_{\bH^r(K)}.
\]
\end{lemma}

\begin{proof}
For $r = 1$, the statement follows from \cite[Lemma~8]{BoffiCDD_06_Dhp}.
For $r = 0$ it is trivial, because
\[
  \inf_{\bq_p \in \bCP^{\rm RT,0}_p(K)} \|\bu - \bq_p\|_{0,K} \le
  \|\bu\|_{0,K}.
\]
Then we obtain the whole range $r \in [0,1]$ via interpolation.
\end{proof}

We will use the above two lemmas to prove the following auxiliary result.

\begin{lemma} \label{lm_RT_aux}
Let $\bv_p \in \bCP^{\rm RT,0}_p(K)$ be such that $\<\bv_p,\bw_p\>_{0,K} = 0$
for any $\bw_p \in \bW_p(K) := \{\bw \in \bCP^{\rm RT,0}_p(K);\; \div\,\bw = 0\}$.
Then there exists $\bv \in \grad\tilde\CH(K)$ such that
\be \label{RT_aux_L^2}
    \|\bv\|_{0,K} \le
    C\, \|\div\,\bv_p\|_{\tH^{-1/2}(K)}
\ee
and
\be \label{RT_aux_H(div)}
    \|\bv - \bv_p\|_{\bH(\div,K)} \le
    C\,p^{\,\eps_0}\,\|\div\,\bv_p\|_{\tH^{-1/2}(K)}
\ee
for any $\eps_0 > 0$ with $C = C(\eps_0) > 0$.
\end{lemma}

\begin{proof}
For given $\bv_p$, we solve the Neumann problem
to find $f \in H^1(K)/{\field{C}}$ such that
\be \label{RT_aux1}
    \< \grad\,f,\grad\,g \>_{0,K} = - \< \div\,\bv_p,g \>_{0,K}\quad
    \forall g \in H^1(K)/{\field{C}}.
\ee
Then we set $\bv := \grad\, f$. One has
\be \label{RT_aux2}
    \div\,\bv = \Delta\,f = \div\,\bv_p \in \CP_{p-1}(K).
\ee
Hence, $f \in \tilde\CH(K)$, $\bv \in \grad\tilde\CH(K)$, and
(\ref{RT_aux_L^2}) holds, because
$\Delta: H^1(K)/{\field{C}} \rightarrow \tH^{-1}(K)$
is an isomorphism.

Note that $\div\,\bv_p \in H^{-1/2+\eps}_{*}(K)$ for any $\eps > 0$.
Therefore, the standard regularity result for problem (\ref{RT_aux1})
reads as (see, e.g., \cite{Grisvard_92_SBV} and cf. Lemma~\ref{lm_LB_reg2}
and Remark~\ref{rm_LB_reg}):
$f \in H^{1+r}(K)$ and
\[
  \|f\|_{H^{1+r}(K)/{\field{C}}} \le C\,\|\div\,\bv_p\|_{H^{-1/2+\eps}(K)}
\]
for any $0 < r \le \frac 12 + \eps$ and for arbitrary $\eps \in (0,1]$.
Then, using the continuity of the gradient as a mapping
$H^{1+r}(K) \rightarrow \bH^r(K)$, we have
\be \label{RT_aux3}
    \|\bv\|_{\bH^r(K)} \le C\,\|f\|_{H^{1+r}(K)/{\field{C}}} \le C\,
    \|\div\,\bv_p\|_{H^{-1/2+\eps}(K)},\quad 0 < r \le 1/2 + \eps,\ \ \eps \in (0,1].
\ee
Since $\bv \in \bH^r(K) \cap \bH_0(\div,K)$, we can
apply the interpolation operator $\Pi^{\div,0}_p$ to define
$\bv^{\div}_p := \Pi^{\div,0}_p\,\bv \in \bCP^{\rm RT,0}_p(K)$.
Recalling that $\Pi^{\div,0}_p$ commutes with the $L^2$-projector
(see (\ref{deRham})) and using (\ref{RT_aux2}), we find that
\[
    \div\,\bv^{\div}_p = \div\,\bv_p = \div\,\bv.
\]
Hence, $(\bv_p - \bv^{\div}_p) \in \bW_p(K)$.
This fact implies the relations
\[
  \<\bv,\bv_p - \bv^{\div}_p\>_{0,K} =
  \<\grad\,f,\bv_p - \bv^{\div}_p\>_{0,K} = 0\quad
  \hbox{and}\quad
  \<\bv_p,\bv_p - \bv^{\div}_p\>_{0,K} = 0,
\]
where the latter equation holds by assumptions on $\bv_p$.
Therefore,
\be \label{RT_aux4}
    \|\bv - \bv_p\|_{0,K} \le
    \|\bv - \bv^{\div}_p\|_{0,K}.
\ee
Since $\bv$ is curl-free, we apply Lemmas~\ref{lm_L^2-stab},~\ref{lm_L^2-approx} and
then use inequality (\ref{RT_aux3}). As a result, we obtain
\bea
    \|\bv - \bv_p^{\div}\|_{0,K}
    & \le &
    C\,p^{1/2}\,\inf_{\bq_p \in \bCP^{\rm RT,0}_p(K)} \|\bv - \bq_p\|_{0,K}
    \nonumber
    \\[5pt]
    & \le &
    C\,p^{1/2-r}\,\|\bv\|_{\bH^r(K)} \le
    C p^{1/2-r} \|\div\,\bv_p\|_{H^{-1/2+\eps}(K)}
    \label{RT_aux5}
\eea
for any $\eps \in (0,1]$.
Then making use of the inverse inequality (see Lemma~\ref{lm_inverse}) we estimate
\be \label{RT_aux6}
     \|\div\,\bv_p\|_{H^{-1/2+\eps}(K)} \le
     C p^{2\eps} \|\div\,\bv_p\|_{H^{-1/2}(K)} \le
     C p^{2\eps} \|\div\,\bv_p\|_{\tH^{-1/2}(K)}.
\ee
Now we again use (\ref{RT_aux2}) and then put together
(\ref{RT_aux4})--(\ref{RT_aux6}). We obtain
\be \label{RT_aux7}
    \|\bv - \bv_p\|_{\bH(\div,K)} = 
    \|\bv - \bv_p\|_{0,K} \le
    C p^{1/2 + 2\eps - r} \|\div\,\bv_p\|_{\tH^{-1/2}(K)}.
\ee
Given an arbitrary $\eps_0 > 0$, we select $\eps = \min\,\{\eps_0,1\}$.
Then (\ref{RT_aux_H(div)}) follows from (\ref{RT_aux7}) by taking $r = \frac 12 + \eps$.
\end{proof}

{\bf $\tilde\bH^{-1/2}(\div)$-conforming interpolation operator.}
Now we proceed to the main goal of this section. Given a vector field
$\bu \in \bH^r(K) \cap \tilde\bH^{-1/2}(\div,K)$ with $r > 0$, we
construct an interpolant $\bu^p = \Pi^{\div,-\frac 12}_p \bu \in \bCP^{\rm RT}_p(K)$.
In particular, $\bu^p$ is defined as the sum of three terms:
\be \label{E^p}
    \bu^p = \bu_1 + \bu^p_2 + \bu^p_3.
\ee
The definition of $\bu_1$ and $\bu^p_2$ follows the construction from
\cite{DemkowiczB_03_pIE}. Let $\bu_1$ be a lowest order interpolant defined as
\be \label{E_1}
    \bu_1 = \sum\limits_{i=1}^{4}
    \Big(\int\limits_{\ell_i} \bu\cdot\bn\,d\sigma\Big)\,\bphi_i,
\ee
where $\bn$ denotes the outward normal unit vector to $K$, and
$\bphi_i$ ($i = 1,\ldots,4$) are the standard basis functions
for $\bCP_1^{\rm RT}(K)$, defined by
\[
  \bphi_i\cdot\bn =
  \cases{
         1 & \hbox{on $\ell_i$},\cr
         0 & \hbox{on $\partial K \backslash \ell_i$}.\cr
        }
\]
For any edge $\ell_i \subset \partial K$ one has
\be \label{E-E_1}
    \int\limits_{\ell_i} (\bu - \bu_1) \cdot \bn\,d\sigma = 0.
\ee
Hence, there exists a function $\psi$, defined on the boundary $\partial K$,
such that
\be \label{psi}
    {\partial\psi\over{\partial \sigma}} =  (\bu - \bu_1) \cdot \bn,\quad
    \psi = 0 \ \ \hbox{at all vertices}.
\ee
Then we define $\psi_2^{\ell_i} \in \CP^0_p(\ell_i)$ by projection
\be \label{psi_2^l}
    \<\psi|_{\ell_i} - \psi_2^{\ell_i}, \phi\>_{\tH^{1/2}(\ell_i)} = 0\quad
    \forall \phi \in \CP^0_p(\ell_i)
\ee
(see Remark~\ref{rem_edge-norms} for the expression of
$\<\cdot,\cdot\>_{\tH^{1/2}(\ell_i)}$).
Extending $\psi_2^{\ell_i}$ by zero from $\ell_i$ onto $\partial K$ (and keeping
its notation), we denote by $\psi_{2,p}^{\ell_i} \in \CP_p(K)$
a polynomial extension of $\psi_2^{\ell_i}$ from $\partial K$ onto $K$, i.e.,
\be \label{psi_2,p^l}
    \psi_{2,p}^{\ell_i} \in \CP_p(K),\quad
    \psi_{2,p}^{\ell_i}|_{\ell_i} = \psi_{2}^{\ell_i},\quad
    \psi_{2,p}^{\ell_i}|_{\partial K\backslash\ell_i} = 0.
\ee
Then we set
\be \label{E_2^p}
    \bu_2^p = \sum\limits_{i=1}^{4} \bu^p_{2,\ell_i},\ \ 
    \hbox{where \ } \bu^p_{2,\ell_i} = \bcurl\, \psi_{2,p}^{\ell_i}.
\ee
The interior interpolant $\bu^p_3$ is a vector bubble function living
in $\bCP^{\rm RT,0}_p(K)$ and satisfying the following system of equations:
\bea
     \lefteqn{
     \<\div(\bu - (\bu_1 + \bu_2^p + \bu_3^p)),\div\,\bv\>_{-1/2,K}
     }
     \nonumber
     \\[3pt]
     &\qquad\qquad = &
     \<\div(\bu - (\bu_1 + \bu_3^p)),\div\,\bv\>_{-1/2,K} = 0\quad
     \forall \bv \in \bCP^{\rm RT,0}_p(K),
     \label{E_3^p_1}
     \\[7pt]
     &   &
     \<\bu - (\bu_1 + \bu_2^p + \bu_3^p),\bcurl\,\phi\>_{0,K} = 0\quad
     \forall \phi \in \CP^{0}_p(K).
     \label{E_3^p_2}
\eea
We note that any polynomial extension satisfying (\ref{psi_2,p^l}) can be
used for the construction of the edge interpolant $\bu_2^p$.
Nevertheless, the interpolant $\bu^p = \bu_1 + \bu^p_2 + \bu^p_3$
is uniquely defined, which follows from (\ref{E_2^p})--(\ref{E_3^p_2}).
It is also easy to see that $\Pi_p^{\div,-\frac 12}$
preserves polynomial vector fields, i.e., $\Pi_p^{\div,-\frac 12} \bv_p = \bv_p$
for any $\bv_p \in \bCP^{\rm RT}_p(K)$.

\begin{remark} \label{rem_diff}
In contrast to the $L^2$-inner product employed to define the
$\bH(\curl)$-conforming interpolation operator (and its
$\bH(\div)$-conforming counterpart) in {\rm \cite{DemkowiczB_03_pIE}},
we use the $\tH^{-1/2}$-inner product in {\rm (\ref{E_3^p_1})}.
\end{remark}

\begin{prop} \label{pr_stab}
For $r > 0$ the operator
\[
  \Pi_p^{\div,-\frac 12}:\;
  \bH^r(K) \cap \tilde\bH^{-1/2}(\div,K) \rightarrow
  \bL^2(K) \cap \tilde\bH^{-1/2}(\div,K)
\]
is well defined and bounded. For arbitrarily small $\eps_0 > 0$
the norm of this operator satisfies
\[
  \Big\|\Pi_p^{\div,-\frac 12}\Big\|_{\CL} \le C\,p^{\,\eps_0},
\]
where $C>0$ is independent of $p$ but depends on $\eps_0$ and $r$,
and $\|\cdot\|_{\CL}$ denotes the operator norm in the space
$\CL\Big(\bH^r(K) \cap \tilde\bH^{-1/2}(\div,K),
 \bL^2(K) \cap \tilde\bH^{-1/2}(\div,K)\Big)$.
\end{prop}

\begin{proof}
Let $\bu \in \bH^r(K) \cap \tilde\bH^{-1/2}(\div,K)$, $r > 0$.
We will study each term on the right-hand side of (\ref{E^p}).
Throughout the proof we denote by $s$ a small parameter such that
$0 < s < \min\,\{\frac 12,r\}$ for given $r>0$.

{\bf Step 1.} Fixing an edge $\ell_i \subset \partial K$ and using a function
\[
  \phi_{i} \in H^{1-s}(K),\quad
  \phi_{i} = \cases{
                         1 & \hbox{on $\ell_i$},\cr
                         0 & \hbox{on $\partial K \backslash \ell_i$}\cr
                        }
\]
as a test function, we integrate by parts to obtain
\beas
     \int\limits_{\ell_i} \bu \cdot \bn\,d\sigma
     & = &
     \int\limits_{\partial K} (\bu \cdot \bn)\,\phi_i\,d\sigma =
     \int\limits_{K} (\div\,\bu)\,\phi_i\,dx +
     \int\limits_{K} \bu \cdot \grad\phi_i\,dx
     \\[3pt]
     & \le &
     \|\div\,\bu\|_{\tH^{-1+s}(K)}\,\|\phi_i\|_{H^{1-s}(K)} +
     \|\bu\|_{\bH^s(K)}\,\|\grad\phi_i\|_{\bH^{-s}(K)}
     \\[3pt]
     & \le &
     C(\phi_i,s)
     \left(
           \|\bu\|_{\bH^r(K)} + \|\div\,\bu\|_{\tH^{-1/2}(K)}
     \right).
\eeas
Note that if $\div\,\bu \in H^{-1+s}(K)$ then an extension to
$\div\,\bu \in \tH^{-1+s}(K)$ exists but is not unique. However,
by assumption $\div\,\bu \in \tH^{-1/2}(K) \subset \tH^{-1+s}(K)$,
which is a unique extension (see~\cite{Mikhailov_08_ATE} for details).
Thus, $\bu_1$ in (\ref{E_1}) is well defined.
Moreover, since $\bu_1$ is a lowest order interpolant, we find
by the equivalence of norms in finite-dimensional spaces that
\be \label{stab1}
    \|\bu_1\|_{\bH(\div,K)} \le
    C \sum\limits_{i=1}^{4} \Big| \int\limits_{\ell_i} \bu \cdot \bn\,d\sigma\Big| \le
    C\left(\|\bu\|_{\bH^r(K)} + \|\div\,\bu\|_{\tH^{-1/2}(K)}\right).
\ee
Let use denote by $\g_{\rm tr}^{-1}$ a right inverse of $\g_{\rm tr}$
with $\g_{\rm tr}^{-1}: H^{1/2-s}(\partial K) \rightarrow H^{1-s}(K)$.
Taking an arbitrary $v \in H^{1/2-s}(\partial K)$ we integrate by parts
similarly as above to estimate
\beas
     \lefteqn{
     \int\limits_{\partial K} (\bu - \bu_1) \cdot \bn\;v\,d\sigma}
     \\
     & \qquad\qquad\le &
     \|\div (\bu - \bu_1)\|_{\tH^{-1+s}(K)}\,\|\g_{\rm tr}^{-1} v\|_{H^{1-s}(K)} +
     \|\bu - \bu_1\|_{\bH^s(K)}\,\|\grad (\g_{\rm tr}^{-1} v)\|_{\bH^{-s}(K)}
     \\[3pt]
     & \qquad\qquad\le &
     C \left(
             \|\bu - \bu_1\|_{\bH^s(K)} + \|\div (\bu - \bu_1)\|_{\tH^{-1+s}(K)}
      \right)
     \|v\|_{H^{1/2-s}(\partial K)}.
\eeas
Hence, $(\bu - \bu_1) \cdot \bn \in H^{-1/2+s}(\partial K)$ and, due to the
finite dimensionality of $\bu_1$, we obtain by using estimate (\ref{stab1}):
\bea
     \|(\bu - \bu_1) \cdot \bn\|_{H^{-1/2+s}(\partial K)}
     & = &
     \sup_{0\not=v\in H^{1/2-s}(\partial K)}
     {|\int_{\partial K} (\bu - \bu_1) \cdot \bn\; v\, d\sigma|
     \over{\|v\|_{H^{1/2-s}(\partial K)}}}
     \nonumber
     \\[3pt]
     & \le &
     C \left(
             \|\bu - \bu_1\|_{\bH^s(K)} + \|\div (\bu - \bu_1)\|_{\tH^{-1+s}(K)}
      \right)
     \nonumber
     \\[3pt]
     & \le &
     C \left(
             \|\bu\|_{\bH^s(K)} + \|\div\,\bu\|_{\tH^{-1+s}(K)} +
             \|\bu_1\|_{\bH(\div,K)}
      \right)
     \nonumber
     \\[3pt]
     & \le &
     C \left(
             \|\bu\|_{\bH^r(K)} + \|\div\,\bu\|_{\tH^{-1/2}(K)}
      \right).
     \label{stab2}
\eea

{\bf Step 2.} From the construction of $\bu_1$ and from the result of
Step~1 we conclude that
\[
  (\bu - \bu_1) \cdot \bn \in H^{-1/2+s}(\partial K),\quad
  \int\limits_{\partial K} (\bu - \bu_1) \cdot \bn\,d\sigma = 0.
\]
Therefore, due to the isomorphism
$\frac{\partial}{\partial \sigma}: H^{1/2+s}(\partial K)/{\field{C}} \rightarrow
 H^{-1/2+s}_{*}(\partial K)$ (see \cite[Lemma~2]{DemkowiczB_03_pIE}),
the function $\psi$ in (\ref{psi}) is well defined, $\psi \in H^{1/2+s}(\partial K)$,
$\psi|_{\ell_i} \in \tH^{1/2}(\ell_i)$ for any edge $\ell_i \subset \partial K$,
and
\be \label{stab3}
    \sum\limits_{i=1}^{4} \|\psi|_{\ell_i}\|_{\tH^{1/2}(\ell_i)} \le
    C\,\sum\limits_{i=1}^{4} \|\psi|_{\ell_i}\|_{H^{1/2+s}_0(\ell_i)} \le
    C\,\|\psi\|_{H^{1/2+s}(\partial K)} \le
    C\,\|(\bu - \bu_1) \cdot \bn\|_{H^{-1/2+s}(\partial K)}.
\ee
Hence, (\ref{psi_2^l}) is uniquely solvable and
\be \label{stab4}
    \|\psi_2^{\ell_i}\|_{\tH^{1/2}(\ell_i)} \le
    C\, \|\psi|_{\ell_i}\|_{\tH^{1/2}(\ell_i)}.
\ee
Furthermore, applying the polynomial extension result
of Babu{\v s}ka and Suri \cite{BabuskaS_87_hpF}, we find the desired
polynomial $\psi_{2,p}^{\ell_i} \in \CP_p(K)$ (see (\ref{psi_2,p^l}))
satisfying
\be \label{stab5}
    \|\psi_{2,p}^{\ell_i}\|_{H^{1}(K)} \le
    C\, \|\psi_2^{\ell_i}\|_{\tH^{1/2}(\ell_i)}.
\ee
Thus, $\bu_2^p$ in (\ref{E_2^p}) is well defined.
Putting together (\ref{stab3})--(\ref{stab5}) we find
\beas
     \|\bu_2^p\|_{0,K}
     & \le &
     C \sum\limits_{i=1}^{4} \|\bcurl\,\psi_{2,p}^{\ell_i}\|_{0,K} \le
     C \sum\limits_{i=1}^{4} \|\psi_{2,p}^{\ell_i}\|_{H^1(K)}
     \\[3pt]
     & \le &
     C \sum\limits_{i=1}^{4} \|\psi_{2}^{\ell_i}\|_{\tH^{1/2}(\ell_i)} \le
     C \sum\limits_{i=1}^{4} \|\psi|_{\ell_i}\|_{\tH^{1/2}(\ell_i)} \le
     C\,\|(\bu - \bu_1) \cdot \bn\|_{H^{-1/2+s}(\partial K)}.
\eeas
Hence, making use of (\ref{stab2}), we obtain
\be \label{stab6}
    \|\bu_2^p\|_{\bH(\div,K)} = \|\bu_2^p\|_{0,K} \le
    C \left(
            \|\bu\|_{\bH^r(K)} + \|\div\,\bu\|_{\tH^{-1/2}(K)}
     \right).
\ee
{\bf Step 3.} The vector bubble function $\bu_3^p$ is uniquely defined by
(\ref{E_3^p_1})--(\ref{E_3^p_2}). To estimate the norms of $\bu^3_p$ and
$\div\,\bu^3_p$ we use the discrete Helmholtz decomposition (\ref{d-dec})
restricted to $\bX_p(K) \equiv \bCP^{\rm RT,0}_p(K)$
\be \label{stab6_1}
    \bu_3^p = \bv_p + \bcurl\,\phi_p,
\ee
where $\phi_p \in \CP_p^0(K)$ and $\bv_p \in \bCP^{\rm RT,0}_p(K)$ is such that
$\<\bv_p,\bw_p\>_{0,K} = 0$ for all $\bw_p \in \bW_p(K)$ (see Lemma~\ref{lm_RT_aux}
for the definition of $\bW_p(K)$).

From (\ref{E_3^p_1}) one has by using the result of Step~1
\bea \label{stab7}
     \|\div\,\bu_3^p\|_{\tH^{-1/2}(K)}
     & \le &
     C\,\|\div(\bu - \bu_1)\|_{\tH^{-1/2}(K)} \le
     C\left(\|\div\,\bu\|_{\tH^{-1/2}(K)} + |\div\,\bu_1|\right)
     \nonumber
     \\[3pt]
     & \le &
     C\left(\|\bu\|_{\bH^{r}(K)} + \|\div\,\bu\|_{\tH^{-1/2}(K)}\right).
\eea
Then we apply Lemma~\ref{lm_RT_aux}:
there exists $\bv \in \grad\tilde\CH(K)$ and an arbitrarily small $\eps_0 > 0$
such that
\[
  \|\bv_p\|_{0,K} \le
  \|\bv - \bv_p\|_{0,K} + \|\bv\|_{0,K} \le
  C\,p^{\,\eps_0}\,\|\div\,\bv_p\|_{\tH^{-1/2}(K)}.
\]
Hence, recalling that $\div\,\bv_p = \div\,\bu^p_3$, we obtain
\be \label{stab8}
    \|\bv_p\|_{0,K} \le
    C\,p^{\,\eps_0}\,\|\div\,\bu_3^p\|_{\tH^{-1/2}(K)}.
\ee
Since $\<\bv_p,\bcurl\,\phi_p\>_{0,K} = 0$, we estimate the norm of
$\bcurl\,\phi_p$ by using (\ref{E_3^p_2}) and by employing the results
of the first two steps:
\be \label{stab9}
    \|\bcurl\,\phi_p\|_{0,K} \le
    \|\bu - \bu_1 - \bu_2^p\|_{0,K} \le
    C\left(\|\bu\|_{\bH^{r}(K)} + \|\div\,\bu\|_{\tH^{-1/2}(K)}\right).
\ee
Combining (\ref{stab7})--(\ref{stab9}) and applying the triangle inequality
we obtain by making use of decomposition (\ref{stab6_1})
\[
  \|\bu_3^p\|_{0,K} + \|\div\,\bu_3^p\|_{\tH^{-1/2}(K)} \le
  C\,p^{\,\eps_0}\left(\|\bu\|_{\bH^{r}(K)} + \|\div\,\bu\|_{\tH^{-1/2}(K)}\right).
\]
To finish the proof it remains to combine the results of the three individual
steps and to apply the triangle inequality to decomposition (\ref{E^p}).
\end{proof}

Proposition~\ref{pr_stab} proves a quasi-stability of the interpolation operator
$\Pi^{\div,-\frac 12}_p$. The following proposition states the commuting diagram
property for $\Pi^{\div,-\frac 12}_p$.

\begin{prop} \label{pr_commute}
For any $\bu \in \bH^r(K) \cap \tilde\bH^{-1/2}(\div,K)$, $r > 0$, there holds
\be \label{P^div_commute}
    \div\,\Big(\Pi^{\div,-\frac 12}_p \bu\Big) = \Pi^{-1/2}_{p-1}(\div\,\bu).
\ee
\end{prop}

\begin{proof}
For any $\varphi \in \CP_{p-1}(K)$ there exists $\bv_p \in \bCP^{\rm RT}_p(K)$
such that $\div\,\bv_p = \varphi$. Therefore, decomposing $\Pi^{\div,-\frac 12}_p \bu$
as in (\ref{E^p}), we need to show that for all $\bv_p \in \bCP^{\rm RT}_p(K)$
there holds
\be \label{commute1}
    \Big\<\div\Big(\bu - \Pi^{\div,-\frac 12}_p \bu\Big), \div\,\bv_p\Big\>_{-1/2,K} =
    \<\div(\bu - (\bu_1 + \bu_3^p)), \div\,\bv_p\>_{-1/2,K} = 0.
\ee
Let us also decompose $\bv_p = \Pi^{\div,-\frac 12}_p \bv_p \in \bCP_p^{\rm RT}(K)$
as in (\ref{E^p}):
\[
  \bv_p = \bv_1 + \bv_2^p + \bv_3^p,\quad
  \div\,\bv_1 = \hbox{const},\quad \div\,\bv_2^p = 0,\quad
  \bv_3^p \in \bCP^{\rm RT,0}_p(K).
\]
Then, recalling (\ref{E_3^p_1}), applying Lemma~\ref{lm_tilde_-1/2_ip},
and integrating by parts, we prove (\ref{commute1}):
\beas
     \lefteqn{
     \Big\<\div\Big(\bu - \Pi^{\div,-\frac 12}_p \bu\Big), \div\,\bv_p\Big\>_{-1/2,K}
     }
     \\
     &\qquad\qquad = &
     \<\div(\bu - \bu_1 - \bu_3^p), \hbox{const}\>_{-1/2,K} +
     \Big\<\div\Big(\bu - \Pi^{\div,-\frac 12}_p \bu\Big), \div\,\bv_3^p\Big\>_{-1/2,K}
     \\
     &\qquad\qquad = &
     \<\div(\bu - \bu_1 - \bu_3^p), \hbox{const}\>_{0,K} =
     \hbox{const}\, \int\limits_{\partial K} (\bu - \bu_1 - \bu_3^p) \cdot \bn\,d\sigma = 0.
\eeas
For the last step we used the fact that $\bu_3^p \cdot \bn|_{\partial K} = 0$
and then applied (\ref{E-E_1}).
\end{proof}

\section{Proof of Theorem~\ref{thm_solve}} \label{sec_proof_solve}
\setcounter{equation}{0}

In this section we prove Theorem~\ref{thm_solve} relying on the abstract
convergence result of Theorem~\ref{thm_abstract}.
Since the discrete decomposition (\ref{d-dec}) is stable with respect
to complex conjugation, one needs to check that assumptions (A1) and (A2) are satisfied.
First, we note that the family $\{\bX_{hp}\}$ of RT-spaces is dense in $\bX^0$.
Since the injection $\bX^0 \subset \bX$ is dense as well
(see, e.g., \cite[Lemma~2.4]{BuffaH_04_CCF}), we conclude that the family
$\{\bX_{hp}\}$ satisfies assumption (A1) of Theorem~\ref{thm_abstract}.

It was mentioned in Section~\ref{sec_dec} that $\bW_{hp}{\subset}\bW$
by construction. Thus, it remains to prove that the subspace
$\bV_{hp}$ defined by (\ref{V_p}) satisfies assumption (\ref{A2}).
In particular, we will show below that there exists a sequence $\{\delta_{hp}\}$,
$\delta_{hp} \rightarrow 0$, such that for any given
$\bv_{hp} \in \bV_{hp}$ there exists $\bv \in \bV$ satisfying
\be \label{Vhp_0}
    \|\bv_{hp} - \bv\|_{\bX} \le \delta_{hp} \|\bv_{hp}\|_{\bX}.
\ee
We prove (\ref{Vhp_0}) for a closed (resp., open) surface $\G$.
The proof consists of 5 steps.

{\bf Step 1: Construction of $\bv$.}
Given $\bv_{hp} \in \bV_{hp}$, we solve the following problem
to find $f \in H^1(\G)/{\field{C}}$ such that
\be \label{Vhp_1}
    \< \gradg\,f,\gradg\,g \> = - \< \divg\,\bv_{hp},g \>\quad
    \forall g \in H^1(\G)/{\field{C}}.
\ee
We set $\bv := \gradg\, f$. Then
\be \label{Vhp_2}
    \divg\,\bv = \divg\,\bv_{hp} \in L^2(\G),
\ee
and, due to Theorem~\ref{thm_H-dec}, there holds $\bv \in \bV$.
Note that $\divg\,\bv_{hp} \in H^{-1/2}(\G)$ (resp.,
$\divg\,\bv_{hp} \in \tH^{-1/2}(\G)$) and $\<\divg\,\bv_{hp},1\> = 0$.
Therefore, the regularity result for problem (\ref{Vhp_1}) reads as
(cf. Lemma~\ref{lm_LB_reg1} (resp., Lemma~\ref{lm_LB_reg2})):
$f \in H^{1+r}(\G)$ with $r = \min\,\{s^{*},\frac 12\} - \eps_1$,
where $s^* > 0$ depends on the geometry of $\G$, $\eps_1 > 0$ is arbitrarily
small. Moreover, using the continuity of $\gradg$, we conclude that
$\bv \in \bH^r_{\;-}(\G)$ and
\be \label{Vhp_3}
    \|\bv\|_{\bH^r_{\;-}(\G)} \le C\,\|f\|_{H^{1+r}(\G)/{\field{C}}} \le C\,
    \|\divg\,\bv_{hp}\|_{\tH^{-1/2}(\G)},\quad
    r = \min\,\{s^{*},\hbox{$\frac 12$}\} - \eps_1
\ee
(here and below we use the convention $\tH^{-1/2}(\G) = H^{-1/2}(\G)$
if $\G$ is closed).

{\bf Step 2: Bounding $\|\bv_{hp} - \bv\|_{\bX}$ by $\|\bv - \bv_{hp}\|_{0,\G}$.}
In view of (\ref{Vhp_2}) we can estimate the norm on the left-hand
side of (\ref{Vhp_0}) as
\bea \label{Vhp_4}
     \|\bv_{hp} - \bv\|_{\bX}
     & = &
     \|\bv_{hp} - \bv\|_{\tilde\bH^{-1/2}_{\|}(\G)} \le 
     C\,\|\bv_{hp} - \bv\|_{\bH^{-1/2+\eps_1}_{\;-}(\G)}
     \nonumber
     \\
     & = &
     C\,\sup_{\bw \in \bH^{1/2-\eps_1}_{\;-}(\G)\setminus\{\bzero\}}
     {\<\bv - \bv_{hp},\bw\>\over{\|\bw\|_{\bH^{1/2-\eps_1}_{\;-}(\G)}}}
\eea
with the same $\eps_1 > 0$ as in (\ref{Vhp_3}).
Let $\bw \in \bH^{1/2-\eps_1}_{\;-}(\G)$. Then, by \cite[Theorem~3.4]{BuffaC_01_TII}
(resp., \cite[Theorem~6.1]{BuffaC_01_TII}), there exists a unique pair
$w_1,\,w_2 \in H^1(\G)/{\field{C}}$
(resp., $w_1 \in H^1(\G)/{\field{C}},\ w_2 \in H^1_0(\G)$) such that
\be \label{Vhp_5}
    \bw = \gradg\, w_1 + \bcurlg\, w_2.
\ee
Moreover, $w_2$ is the solution of the problem
\[
  - \Delta_\G\,w_2 = \curlg\,\bcurlg\, w_2 = \curlg\,\bw.
\]
Since $\curlg\,\bw \in H^{-1/2-\eps_1}(\G)$ and $\<\curlg\,\bw,1\> = 0$
(resp., $\curlg\,\bw \in H^{-1/2-\eps_1}(\G)$),
we apply  Lemma~\ref{lm_LB_reg1} (resp., the Dirichlet analog of Lemma~\ref{lm_LB_reg2})
to prove that $w_2 \in H^{1+r}(\G)$ (with the same $r$ as in (\ref{Vhp_3}))
and there holds
\be \label{Vhp_6}
    \ba{c}
    \|w_2\|_{H^{1+r}(\G)/{\field{C}}} \le C\,
    \|\curlg\,\bw\|_{H^{-1/2-\eps_1}(\G)} \le
    C\,\|\bw\|_{\bH^{1/2-\eps_1}_{\;-}(\G)}
    \\[11pt]
    (\hbox{resp.,}\quad
    \|w_2\|_{H^{1+r}(\G)} \le
    C\,\|\bw\|_{\bH^{1/2-\eps_1}_{\;-}(\G)}).
    \ea
\ee
Therefore, $w_2 \in H^{1+r}(\G_j)$, $r>0$, for any element $\G_j$, and one can
apply the $H^1$-conforming interpolation operator $\Pi_p^1$
to find a continuous piecewise polynomial $w_{2}^{hp}$ defined on $\G$ such that
$\hat w_{2,j}^{hp} := w_{2}^{hp}|_{\G_j} \circ T_j = \Pi_p^1 \hat w_{2,j} \in \CP_p(K)$
(here $\hat w_{2,j} := w_{2}|_{\G_j} \circ T_j$). To estimate the $H^1$-semi-norm of the
error $(w_2 - w_2^{hp})$, we apply the standard scaling argument
(cf. \cite[Theorem~4.3.2]{Ciarlet_78_FEM}) and the approximation result for $\Pi_p^1$
on the reference element (see \cite[Theorem~2]{DemkowiczB_03_pIE}):
\bea
    |w_2 - w_{2}^{hp}|_{H^1(\G_j)}
    & \le &
    C\,\Big\|\hat w_{2,j} - \Pi_p^1 \hat w_{2,j}\Big\|_{H^1(K)}
    \nonumber
    \\[5pt]
    & = &
    C\,\Big\|\hat w_{2,j} - \hat\varphi_p -
             \Pi_p^1 (\hat w_{2,j} - \hat\varphi_p)\Big\|_{H^1(K)}\qquad
    (\forall \hat\varphi_p \in \CP_p(K))
    \nonumber
    \\[5pt]
    & \le &
    C\,p^{-(r-\eps_2)}\,
    \inf_{\hat\varphi_p \in \CP_p(K)} \|\hat w_{2,j} - \hat\varphi_p\|_{H^{1+r}(K)},\quad
    0 < \eps_2 < r < 1/2.\qquad
    \label{Vhp_6.5}
\eea
Let $s = 1,2$. Using Theorem~3.1.1 of \cite{Ciarlet_78_FEM} and the scaling argument,
one has
\beas
     \inf_{\hat\varphi_p \in \CP_p(K)} \|\hat w_{2,j} - \hat\varphi_p\|_{H^{s}(K)}
     & \le &
     \inf_{\hat\varphi \in \CP_{s-1}(K)} \|\hat w_{2,j} - \hat\varphi\|_{H^{s}(K)}
     \\[5pt]
     & \le &
     C\,|\hat w_{2,j}|_{H^{s}(K)} \le
     C\,h^{s-1}\,\|w_2\|_{H^{s}(\G_j)}.
\eeas
Therefore, by interpolation,
\[
  \inf_{\hat\varphi_p \in \CP_p(K)} \|\hat w_{2,j} - \hat\varphi_p\|_{H^{1+r}(K)} \le
  C\,h^{r}\,\|w_2\|_{H^{1+r}(\G_j)},
\]
and from (\ref{Vhp_6.5}) we conclude that
\[
  |w_2 - w_{2}^{hp}|_{H^1(\G_j)} \le
  C\,h^{r}\,p^{-(r-\eps_2)}\,\|w_2\|_{H^{1+r}(\G_j)}.
\]
Hence, for a closed (resp., open) surface $\G$ we obtain
\be \label{Vhp_7}
    \ba{c}
    |w_2 - w_2^{hp}|_{H^1(\G)} =
    \bigg(\sum\limits_{j} |w_2 - w_2^{hp}|^2_{H^1(\G_j)}\bigg)^{1/2} \le
    C\, h^r\, p^{-(r-\eps_2)}\, \|w_2\|_{H^{1+r}(\G)/{\field{C}}}
    \\[11pt]
    (\hbox{resp.,}\quad
    |w_2 - w_2^{hp}|_{H^1(\G)} \le
    C\, h^r\, p^{-(r-\eps_2)}\, \|w_2\|_{H^{1+r}(\G)}).
    \ea
\ee
In addition, for the case of an open surface, we note that $w_2^{hp}$ vanishes on $\partial\G$.
Then recalling the commuting diagram property for $\Pi^1_p$
(see (\ref{deRham})) and the definition of $\bX_{hp}$ (see (\ref{Xp})),
we conclude that $\bcurlg\,w_2^{hp} \in \bX_{hp} \subset \bX$.
Moreover, $\bcurlg\,w_2^{hp} \in \bW_{hp} \subset (\bW \cap \bL^2_t(\G))$.
This fact together with the $\bL^2_t(\G)$-orthogonalities
$\bV \perp (\bW \cap \bL^2_t(\G))$ and $\bV_{hp} \perp \bW_{hp}$ implies
\be \label{Vhp_8}
    \<\bv - \bv_{hp},\bcurlg\,w_2^{hp}\> = 0.
\ee
Now we use (\ref{Vhp_2}), (\ref{Vhp_5}), (\ref{Vhp_6}), (\ref{Vhp_7}), (\ref{Vhp_8})
(and (\ref{X_property}), if $\G$ is an open surface)
to prove for any $\bw \in \bH^{1/2-\eps_1}_{\;-}(\G)$
\beas
     \<\bv - \bv_{hp},\bw\>
     & = &
     \<\bv - \bv_{hp},\gradg\,w_1\> + \<\bv - \bv_{hp},\bcurlg\,w_2\>
     \\[3pt]
     & = &
     -\<\divg(\bv - \bv_{hp}),w_1\> + \<\bv - \bv_{hp},\bcurlg(w_2 - w_2^{hp})\>
     \\[3pt]
     & \le &
     \|\bv - \bv_{hp}\|_{0,\G}\,|w_2 - w_2^{hp}|_{H^1(\G)} \le
     C\,h^r\,p^{-(r-\eps_2)}\,\|\bv - \bv_{hp}\|_{0,\G}\,\|\bw\|_{\bH^{1/2-\eps_1}_{\;-}(\G)}.
\eeas
Using this estimate in (\ref{Vhp_4}) we find
\be \label{Vhp_9}
    \|\bv_{hp} - \bv\|_{\bX} \le 
    C\,h^r\,p^{-(r-\eps_2)}\,\|\bv - \bv_{hp}\|_{0,\G},
\ee
where $r = \min\,\{s^{*},\frac 12\}-\eps_1$, \ $\eps_1 > 0$, \ $\eps_2 \in (0,r)$ \
($\eps_1,\, \eps_2$ can be arbitrarily small).

{\bf Step 3: Bounding $\|\bv - \bv_{hp}\|_{0,\G}$ by $\|\bv - \Pi^{\div,-1/2}_p\bv\|_{0,\G}$.}
We recall that $\bv \in \bH^r_{\;-}(\G) \cap \bH(\divg,\G)$,
$r \in (0,\frac 12)$. Therefore, $\bv|_{\G_j} \in \bH^r(\G_j) \cap \bH(\divg,\G_j)$
for each element $\G_j$ and we can define $\bv_{hp}^{\div} \in \bX_{hp}$ such that
\[
  \CM^{-1}_{j}\Big(\bv_{hp}^{\div}|_{\G_j}\Big) = 
  \Pi^{\div, -\frac 12}_p\Big(\CM^{-1}_{j}(\bv|_{\G_j})\Big),
\]
where $\CM_j$ is the Piola transform (see (\ref{Piola})).

Using commutativity (\ref{P^div_commute}) and the fact that
$\divg\,\bv = \divg\,\bv_{hp}$ is a piecewise polynomial on $\G$,
we find
\[
    \divg\,\bv^{\div}_{hp} = \divg\,\bv = \divg\,\bv_{hp}.
\]
Hence, $(\bv_{hp} - \bv^{\div}_{hp}) \in \bW_{hp} \subset (\bW \cap \bL^2_t(\G))$,
and using again the orthogonalities $\bV \perp (\bW \cap \bL^2_t(\G))$,
$\bV_{hp} \perp \bW_{hp}$ we have
\[
  \<\bv - \bv_{hp},\bv_{hp} - \bv^{\div}_{hp}\> = 0.
\]
Therefore,
\be \label{Vhp_10}
    \|\bv - \bv_{hp}\|_{0,\G} \le \|\bv - \bv^{\div}_{hp}\|_{0,\G}.
\ee

{\bf Step 4: Estimating $\|\bv - \bv^{\div}_{hp}\|_{0,\G}$.}
Using (\ref{Piola_2}) and the quasi-stability of $\Pi^{\div,-\frac 12}_p$
(see Proposition~\ref{pr_stab}), we estimate
\beas
     \|\bv^{\div}_{hp}|_{\G_j}\|_{0,\G_j}
     & = &
     \Big\|\CM_j\Big(\Pi^{\div,-\frac 12}_p(\CM_j^{-1}(\bv|_{\G_j}))\Big)\Big\|_{0,\G_j} \le
     C\,\Big\|\Pi^{\div,-\frac 12}_p(\CM_j^{-1}(\bv|_{\G_j}))\Big\|_{0,K}
     \\[5pt]
     & \le &
     C\,p^{\eps_0}\,
     \Big(
     \|\CM_j^{-1}(\bv|_{\G_j})\|_{\bH^r(K)} +
     \|\div(\CM_j^{-1}(\bv|_{\G_j}))\|_{\tH^{-1/2}(K)}
     \Big)
     \\[5pt]
     & \le &
     C\,p^{\eps_0}\,
     \Big(
     \|\CM_j^{-1}(\bv|_{\G_j})\|_{\bH^r(K)} +
     \|\div(\CM_j^{-1}(\bv|_{\G_j}))\|_{H^{-1/2+\eps_3}(K)}
     \Big),
\eeas
where $\eps_0,\eps_3 > 0$ are arbitrarily small.
Hence, recalling that $\div(\CM_j^{-1}(\bv|_{\G_j})) = \div(\CM_j^{-1}(\bv_{hp}|_{\G_j}))$
is a polynomial on $K$, we make use of the inverse inequality (see Lemma~\ref{lm_inverse})
and then apply (\ref{Piola_3}), (\ref{Piola_4}):
\beas
     \|\bv^{\div}_{hp}|_{\G_j}\|_{0,\G_j}
     & \le &
     C\,p^{\eps_0+2\eps_3}\,
     \Big(
     \|\CM_j^{-1}(\bv|_{\G_j})\|_{\bH^r(K)} +
     \|\div(\CM_j^{-1}(\bv_{hp}|_{\G_j}))\|_{H^{-1/2}(K)}
     \Big)
     \\[5pt]
     & \le &
     C\,p^{\eps_0+2\eps_3}\,
     \Big(\|\bv\|_{\bH^r(\G_j)} + h^{1/2}\,\|\divg\,\bv_{hp}\|_{H^{-1/2}(\G_j)}\Big).
\eeas
Therefore,
\bea
    \|\bv - \bv^{\div}_{hp}\|_{0,\G}
    & \le &
    \|\bv\|_{0,\G} +
    \bigg(\sum\limits_{j} \|\bv^{\div}_{hp}|_{\G_j}\|^2_{0,\G_j}\bigg)^{1/2}
    \nonumber
    \\[3pt]
    & \le &
    \|\bv\|_{0,\G} +
    C\,p^{\eps_0+2\eps_3}\,
    \bigg(
          \sum\limits_{j}
          \Big(\|\bv\|^2_{\bH^r(\G_j)} + h\,\|\divg\,\bv_{hp}\|^2_{H^{-1/2}(\G_j)}\Big)
    \bigg)^{1/2}\qquad
    \nonumber
    \\[3pt]
    & \le &
    C\,p^{\,\eps_0 + 2\eps_3}
    \left(
    \|\bv\|_{\bH^r_{\;-}(\G)} + h^{1/2}\,\|\divg\,\bv_{hp}\|_{H^{-1/2}(\G)}
    \right).
    \label{Vhp_11}
\eea

{\bf Step 5: Conclusion.}
First, we select all small parameters $\eps_k$ ($k = 0,1,2,3$)
such that $\eps_0+\eps_1+\eps_2+2\eps_3 = \eps$ for any given
$\eps \in (0,\min\,\{s^{*},\frac 12\})$. Then putting together the results of all
previous steps, i.e., estimates (\ref{Vhp_3}) and (\ref{Vhp_9})--(\ref{Vhp_11}),
we obtain
\[
  \|\bv_{hp} - \bv\|_{\bX} \le
  C\left(\hbox{$\frac hp$}\right)^{\min\,\{s^{*},\frac 12\} - \eps}\,\|\divg\,\bv_{hp}\|_{\tH^{-1/2}(\G)}\le
  C\left(\hbox{$\frac hp$}\right)^{\min\,\{s^{*},\frac 12\} - \eps}\,\|\bv_{hp}\|_{\bX}.
\]
This proves (\ref{Vhp_0}). Therefore, the subspace $\bV_{hp}$ satisfies (\ref{A2}).
Thus we have shown that the discrete decomposition (\ref{d-dec})
verifies assumption (A2) of Theorem~\ref{thm_abstract} in the framework of the
$hp$-version of the BEM with quasi-uniform meshes, and the proof is finished.

\section{Concluding remarks} \label{sec_conclusion}
\setcounter{equation}{0}

The main result of this paper -- the convergence of the $hp$-BEM with quasi-uniform meshes
(and thus the convergence of the $h$- and the $p$- versions as particular cases)
for the EFIE -- is proved for a sequence of Raviart-Thomas spaces on quadrilateral
elements, which can be parallelograms, convex quadrilaterals or curvilinear ones.
To that end it was essential to introduce and analyse
a new $\tilde\bH^{-1/2}(\div)$-conforming interpolation operator $\Pi_p^{\div,-\frac 12}$
on the reference element. To show the stability of $\Pi_p^{\div,-\frac 12}$ (with
respect to polynomial degrees) we relied, in particular, on the discrete inf-sup condition
(see Lemma~\ref{lm_inf-sup}), which was established in~\cite{BoffiCDD_06_Dhp}.

{\em The case of triangular elements.}
We note that Lemma~\ref{lm_inf-sup} is not available for the reference triangle.
However, the corresponding result has been conjectured and numerically evidenced
in~\cite{BoffiDC_03_DCp} for $\bH(\curl)$-conforming N{\' e}d{\' e}lec
elements of the second type. In our $\bH(\div)$-conforming settings it means
that the conjectured result of Lemma~\ref{lm_inf-sup} is numerically
justified for Brezzi-Douglas-Marini spaces on the reference triangle
(this is due to the isomorphism of the curl and the $\div$ operator in 2D).
Since all remaining arguments in the construction of $\Pi_p^{\div,-\frac 12}$ and in
the proofs of Propositions~\ref{pr_stab},~\ref{pr_commute} and Theorem~\ref{thm_solve}
are valid for BDM-spaces on triangles (cf. also~\cite{DemkowiczB_03_pIE}),
we conclude that the main result of the paper holds for these spaces, provided
that the discrete inf-sup condition discussed above is true.

{\em Regularity results and error analysis.}
In this paper we do not present the regularity result for the problem under
consideration, neither perform an a priori error analysis for the $hp$-BEM.
However, we note that in~\cite{BespalovH_NpB} we derived explicit formulas for
typical singularities inherent to the solution of the EFIE on piecewise smooth
(open or closed) Lipschitz surfaces. The $p$-approximation analysis of these singularities
(including the least regular ones) was performed in~\cite{BespalovH_NpB} on a plane
open surface. This analysis relied on our results for the Laplacian,
see~\cite{BespalovH_05_pBE,BespalovH_07_pBW,Bespalov_NPA}, by using continuity
properties of the surface (vector) curl operator. Since those results and properties
(in corresponding spaces) are valid for polyhedral and piecewise plane open surfaces,
the proof of an a priori error estimate for the $p$-BEM from {\rm \cite{BespalovH_NpB}}
carries over without essential modifications to the more general case considered
in the present paper. We stress that this $p$-approximation result holds only for
affine families of meshes. The error analysis of the $hp$-version on quasi-uniform
meshes presents more difficulties. In particular, the involved Sobolev spaces of
negative order are not scalable under affine transformations. Therefore, this analysis
is not a trivial extension of our $p$-approximation results.

{\em Non-affine quadrilateral meshes.}
Although the main convergence result of the paper holds for non-affine
quadrilateral meshes, it is not clear if the RT-spaces in this case
would provide optimal approximations in the energy norm for the EFIE.
For instance, in the $h$-version of the finite element method the degradation
of convergence rates (in the $\bH(\div)$-norm) was observed for RT-elements on
general convex quadrilaterals (see~\cite{ArnoldBF_05_QHF}).
An alternative family of finite elements on non-affine meshes
was introduced in \cite{ArnoldBF_05_QHF} and was shown to provide
optimal $h$-approximation order in $\bH(\div)$. These Arnold-Boffi-Falk (ABF) elements
can also be used in the BEM for the EFIE. However, the solvability, convergence, and
a priori error estimation of the $h$-, $p$-, or $hp$-BEM
with ABF-elements are open problems.

\bibliographystyle{siam}
\bibliography{bib,heuer,fem}

\end{document}